\documentclass[11pt]{amsart} 
\usepackage{verbatim, latexsym, amssymb, amsmath,color}
\usepackage{epsfig}

\def\R{\mathbb R}
\def\RP{\mathbb {RP}}
\def\N{\mathbb N}
\def\Z{\mathbb Z}

\theoremstyle{remark}

\theoremstyle{definition}

\title[Existence of  infinitely many minimal hypersurfaces]{Existence of infinitely many minimal hypersurfaces in positive Ricci curvature}
\author{Fernando C. Marques and Andr\'e Neves}
\address{Fine Hall\\Princeton University \\ Princeton, NJ 08544}
\email{coda@math.princeton.edu}
\address{Imperial College London\\ Huxley Building \\ 180 Queen's Gate \\ London SW7 2RH \\ United Kingdom}
\email{a.neves@imperial.ac.uk}
\thanks{The first author was partly supported by CNPq-Brazil, FAPERJ, CNRS and Edital Brazil-France: ANR-11-IS01-0002. The second author was partly supported by Marie Curie IRG Grant and ERC Start Grant.}

\begin{document}
\maketitle
\begin{abstract}
In the early 1980s, S. T. Yau conjectured that any compact  Riemannian three-manifold admits an infinite number of closed immersed  minimal surfaces. We use min-max theory for the area functional to prove this conjecture
in the positive Ricci curvature setting. More precisely, we show that 
every compact Riemannian manifold with positive Ricci curvature and dimension at most seven contains  infinitely many  smooth, closed, embedded   minimal hypersurfaces. 

In the last section we mention some open problems related with the geometry of these minimal hypersurfaces.
\end{abstract}
\section{Introduction}

A foundational question in Differential Geometry, asked by Poincar\'e \cite{poincare}, is whether every closed Riemann surface admits a closed geodesic. 
If the surface is not simply connected then we can minimize length in a nontrivial homotopy class and produce a closed geodesic. Therefore the question becomes considerably more interesting on a sphere, and the first breakthrough was due to Birkhoff \cite{birkhoff}  who used min-max methods to find a closed geodesic for any  metric on a two-sphere. 

Later, in a remarkable work, Lusternik and Schnirelmann \cite{lusternik}  showed that every metric on a $2$-sphere admits three simple (embedded) closed geodesics (see also \cite{ballman, grayson, jost2, klingenberg, lusternik2, taimanov}). This suggests the question of whether we can find  an infinite number of geometrically distinct closed geodesics in any closed surface. It is not hard to find infinitely many closed geodesics
when the genus of the surface is positive.
The case of the sphere was finally settled by Franks \cite{franks} and Bangert \cite{bangert}. Their  works combined  imply that every metric on a two-sphere admits an infinite number of closed geodesics.  Later, Hingston \cite{hingston} estimated the number of closed geodesics of length at most $L$ when $L$ is very large.

Likewise, one can ask whether every closed Riemannian manifold admits a closed minimal hypersurface. Using min-max methods, and building on earlier work of Almgren, Pitts \cite{pitts} proved that  every compact Riemannian $(n+1)$-manifold with $n\leq 5$ contains a smooth, closed, embedded  minimal hypersurface. Later, Schoen and Simon \cite{schoen-simon} extended this result to any dimension, proving the existence of a closed,  embedded minimal hypersurface with a singular set of Hausdorff codimension at least $7$. 

Motivated by these results, Yau conjectured in \cite{yau1} (first problem in the Minimal Surfaces section) that every compact  Riemannian three-manifold admits an infinite number of smooth, closed, immersed minimal surfaces. The main purpose of this paper is to  prove this conjecture in the positive Ricci curvature setting. More generally, we establish the existence of infinitely many  smooth, closed, embedded, minimal hypersurfaces for manifolds that satisfy a Frankel-type property and have dimension less than or equal to 7.

The main result of  this paper is:

\subsection{Main Theorem}\label{main.thm}
{\em Let $(M^{n+1},g)$ be a compact Riemannian manifold, $3 \leq (n+1)\leq 7$. Then either
\begin{itemize}
\item[(i)] there exists a disjoint collection $\{\Sigma_1,\dots, \Sigma_{n+1}\}$ of $(n+1)$ connected closed smooth embedded minimal hypersurfaces, 
\item[(ii)] or there exist infinitely many connected closed smooth embedded minimal hypersurfaces.
\end{itemize}
}

\subsection{Corollary}\label{multiple.existence}
{\em Every compact Riemannian manifold $(M^{n+1},g)$ of dimension $3\leq (n+1)\leq 7$ contains at least $(n+1)$ connected closed smooth embedded minimal hypersurfaces.}

\medskip

\subsection{Definition} We say that a Riemannian manifold $(M,g)$ satisfies the {\it embedded Frankel property} if  any two closed, smooth embedded minimal hypersurfaces of $M$ intersect each other.

\medskip

\subsection{Corollary}\label{frankel.existence}{\em Let $(M,g)$ be a compact Riemannian manifold of dimension $3\leq (n+1)\leq 7$. Suppose that $M$ satisfies the embedded Frankel property. Then $M$ contains an infinite number of distinct closed, smooth embedded, minimal hypersurfaces.
}

\medskip

Since manifolds of positive Ricci curvature satisfy the embedded Frankel property \cite{frankel}, we derive the following corollary:

\subsection{Corollary}\label{ricci.thm}{\em Let $(M,g)$ be a compact  Riemannian $(n+1)$-manifold with $3\leq (n+1)\leq 7$. If the Ricci curvature of $g$ is positive, then $M$ contains an infinite number of distinct  closed, smooth embedded, minimal hypersurfaces.
}

\subsection{Remark:} { The counterparts of the Main Theorem, Corollary \ref{multiple.existence}, Corollary \ref{frankel.existence} and Corollary \ref{ricci.thm} in higher dimensions hold true if one allows the minimal hypersurfaces to be smooth outside sets of codimension 7. These extensions can be found in \cite{marques-neves-cycles}.}

\medskip

The proof of the Main Theorem uses the Almgren-Pitts min-max theory for the volume functional, combined with ideas from Lusternik-Schnirelmann theory. The idea is to  apply min-max theory to the high-parameter families of hypersurfaces (mod $2$ cycles) constructed by Guth in \cite{guth}. We give an informal overview of the proof  at the end of this section.

The Almgren-Pitts min-max theory does not produce closed geodesics when the ambient is a two-dimensional surface ($n=1$). The min-max varieties can be stationary geodesic networks with point singularities, since they satisfy the almost minimizing in annuli condition (\cite{pitts0}). In fact it is well-known that there are ellipsoids in $\mathbb{R}^3$ with exactly three embedded closed geodesics.

 In \cite{kapouleas1, kapouleas2} Kapouleas describes in detail an alternative approach to construct an infinite number of  embedded minimal surfaces in a three-manifold with a generic metric by either desingularizing two intersecting minimal surfaces or by doubling an existing unstable minimal surface. Note that  for $S^3$ with a metric of positive Ricci curvature, White \cite{white2} showed the existence  of two distinct embedded minimal spheres, which must intersect by \cite{frankel} and are necessarily unstable.
 
 The minimal hypersurfaces obtained via our construction  have, conjecturally, area  tending to infinity and thus should  be different from the minimal surfaces proposed by Kapouleas.
 
Rubinstein \cite{rubinstein} outlined an argument to produce  an infinite number of minimal immersed surfaces in any hyperbolic $3$-manifold with finite volume. He assumes, among other things, that minimal surfaces produced  from Heegaard splittings via min-max methods have index one but this remains an open problem.

Some other conditions are known to imply the embedded Frankel property. For instance, any closed Riemannian manifold $(M^{n+1},g)$, $2\leq n \leq 6$, that does not admit compact, embedded minimal hypersurfaces with stable two-sided covering satisfies the embedded Frankel property. This follows from the same argument as in Theorem 9.1 of \cite{meeks-perez-ros}. Hence:

\subsection{Corollary}\label{stable.thm}{\em Let $(M,g)$ be a compact  Riemannian $(n+1)$-manifold with $2\leq n\leq 6.$ Suppose that $(M,g)$ contains no  closed, embedded minimal hypersurfaces with stable two-sided covering. Then $M$ contains an infinite number of distinct smooth, closed, embedded, minimal hypersurfaces. 
}

\medskip

\subsection{Remark:} The families we use in this paper have analogues for the case of compact manifolds with boundary. In fact, these are the families (of relative cycles) considered by Guth \cite{guth} in the unit ball.  Once the Almgren-Pitts theory
is adapted to that setting, the arguments of this paper should lead to the existence of infinitely many distinct smooth, properly embedded, free boundary minimal hypersurfaces, provided the ambient manifold satisfies a Frankel property. We refer the reader to the paper of Li and Zhou \cite{li-zhou} for details. The Frankel property in the free boundary setting is established  in  Lemma 2.4 of \cite{li-fraser} for compact manifolds with nonnegative Ricci curvature and strictly convex boundary. Geodesic balls with a rotationally symmetric metric also satisfy this property. This last fact follows by using ambient rotations and applying   the maximum principle, and has been pointed out to us by Harold Rosenberg.

\subsection{Overview of the proof:} {The homotopy groups of the space of modulo $2$  $n$-cycles in $M$, $\mathcal Z_{n}(M,\Z_2)$, can be computed through the work of Almgren \cite{almgren}. It follows that   all homotopy groups vanish but the first one: $\pi_1(\mathcal Z_{n}(M,\Z_2)) = \Z_2$, just like in   $\RP^{\infty}$. We consider the generator $\bar \lambda \in H^1(\mathcal Z_{n}(M,\Z_2), \Z_2).$  


Guth \cite{guth} and Gromov \cite{gromov0, gromov, gromov2} have studied  continuous maps $\Phi$ from a simplicial complex $X$ into $\mathcal Z_{n}(M,\Z_2)$ that detect $\bar \lambda^p$, in the sense that $\Phi^{*}(\bar \lambda^p)\neq 0$. In particular, it follows from their construction that for every $p\in \N$ there exists a map  $\Phi$ that detects $\bar\lambda^p$ (with $X=\RP^p$) and such that
 $$\sup_{x\in \RP^p}{\bf M}(\Phi(x)) \leq  C p^{\frac {1}{n+1}},$$
 where $C$ depends only on $M$.  Here ${\bf M}(T)$ denotes the mass of $T$.} Guth's construction  was based on an elegant bend--and--cancel argument that we present in  Section \ref{guth} for the reader's convenience.

 Thus, denoting by $\mathcal P_p$ the space of all maps  that detect $\bar\lambda^p$, we have (see also \cite[Appendix 3]{guth})
 \begin{equation}\label{sublinear}
 \omega_p:=\inf_{\Phi\in \mathcal P_p}\sup_{x\in {\rm dmn}(\Phi)} {\bf M}(\Phi(x))\leq Cp^{\frac{1}{n+1}},
 \end{equation}
 where ${\rm dmn}(\Phi)$ stands for the domain of $\Phi$.

 In Section \ref{equality} we use Lusternik-Schnirelmann theory to show that if $\omega_p=\omega_{p+1}$ then there are infinitely many embedded minimal hypersurfaces.  
 
 The main theorem is proven by contradiction, where we assume that there exist only finitely many smooth, closed, embedded minimal hypersurfaces. Then  $\{\omega_p\}_{p\in \N}$ is strictly increasing and,  under the Frankel condition, each min-max volume $\omega_p$ must be  achieved by  a {\em connected}, closed, embedded minimal hypersurface with some integer multiplicity. In Section \ref{proof.main.theorem} we use this to show that $\omega_p$ must grow linearly  in $p$ and this  is in contradiction with the sublinear growth of $\omega_p$ in $p$ given in \eqref{sublinear}.
 
Sections \ref{almgren.pitts.theory}, \ref{almgren.homo}, \ref{minmax.families} are used to set up and state  the results we need from Almgren--Pitts Min-max Theory. The need for a careful and detailed account in these sections comes from the fact that Almgren--Pitts theory uses the mass norm in ${\mathcal Z}_k(M;\Z_2)$ and sequences of discrete maps, while the elements in $\mathcal P_p$ are continuous maps into ${\mathcal Z}_k(M;\Z_2)$ with respect to the flat topology. Thus it is important to have the technical tools that allow us to move  from one concept to the other. 


{{\bf Acknowledgements:}   Part of this work was done during the first author's stay in Paris. He is  grateful to \'{E}cole Polytechnique, \'{E}cole Normale Sup\'{e}rieure and Institut Henri Poincar\'{e} for the hospitality.}

\section{Almgren-Pitts Min-Max Theory}\label{almgren.pitts.theory}


Let $(M,g)$ be an orientable compact Riemannian $(n+1)$-manifold, possibly with boundary $\partial M$. We  assume that $M$ is
isometrically embedded into some Euclidean space  $\mathbb{R}^L$.  

Let $X$ be a  cubical subcomplex of the  $m$-dimensional cube $I^m=[0,1]^m$. Each $k$-cell of $I^m$ is of the form  $\alpha_1 \otimes \cdots \otimes \alpha_m$, where $\alpha_i \in \{0,1,[0,1]\}$ for every $i$ and $\sum \textrm{dim}(\alpha_i)=k$. Notice that every polyhedron is homeomorphic to the support of some cubical subcomplex of this type \cite[Chapter 4]{panov}.   

  We now describe the necessary and obvious modifications to the Almgren-Pitts Min-Max Theory so that the $m$-dimensional cube $I^m$ is replaced by $X$ as the parameter space.
  
  \subsection{Basic notation}\label{gmt}

The spaces we will work with in this paper are:
\begin{itemize}
\item the space ${\bf I}_k(M;\Z_2)$  of $k$-dimensional mod 2 flat chains    in $\mathbb{R}^L$ with support contained  in $M$ (see \cite[4.2.26]{federer} for more details);
\item the space ${\mathcal Z}_k(M;\Z_2)$ (${\mathcal Z}_k(M,\partial M;\Z_2)$) of mod 2 flat chains  $T \in {\bf I}_k(M;\Z_2)$ with  $\partial T=0$ (${\rm spt }(\partial T) \subset \partial M $);
\item the closure $\mathcal{V}_k(M)$, in the weak topology, of the space of $k$-dimensional rectifiable varifolds in $\mathbb{R}^L$ with support contained in $M$. The space of integral rectifiable $k$-dimensional varifolds with support contained in $M$ is denoted by $\mathcal{IV}_k(M)$.
\end{itemize}

Given $T\in {\bf I}_k(M;\Z_2)$,  we denote by $|T|$ and $||T||$ the integral varifold   and the Radon measure in $M$ associated with $T$, respectively;
 given $V\in \mathcal{V}_k(M)$, $||V||$ denotes the Radon measure in $M$ associated with $V$.  
 If $U\subset M$ is an open set of finite perimeter, we abuse notation and denote the associated current in ${\bf I}_{n+1}(M;\Z_2)$  by $U$.
 
 The above spaces come with several relevant metrics. The  {\it flat metric} and the {\it mass} of $T \in {\bf I}_k(M;\Z_2)$, denoted by $\mathcal{F}(T)$ and ${\bf M}(T)$,  are defined in \cite[page 423]{federer} and  \cite[page 426]{federer},  respectively.  The  ${\bf F}$-{\it metric} on $\mathcal{V}_k(M)$ is defined in  {Pitts book} \cite[page 66]{pitts} and   induces the varifold weak topology on $\mathcal{V}_k(M)$.
Finally,  the ${\bf F}$-{\it metric} on ${\bf I}_k(M;\Z_2)$ is defined by
$$ {\bf F}(S,T)=\mathcal{F}(S-T)+{\bf F}(|S|,|T|).$$

We assume that  ${\bf I}_k(M;\Z_2)$, ${\mathcal Z}_k(M;\Z_2)$, and ${\mathcal Z}_k(M,\partial M;\Z_2)$  have the topology induced by the flat metric. When endowed with
the topology of the mass norm, these spaces will be denoted by  ${\bf I}_k(M;{\bf M};\Z_2)$, ${\mathcal Z}_k(M;{\bf M};\Z_2)$, and ${\mathcal Z}_k(M, \partial M;{\bf M};\Z_2)$, respectively. The space $\mathcal{V}_k(M)$ is considered with the weak topology of varifolds.
Given $\mathcal{A,B}\subset \mathcal{V}_k(M)$, we also define
  $${\bf F}(\mathcal{A},\mathcal{B})=\inf\{{\bf F}(V,W):V\in \mathcal{A}, W\in \mathcal{B}\}.$$
  
   For each $j\in \N$, $I(1,j)$ denotes the cube complex on $I^1$  whose $1$-cells and $0$-cells (those are sometimes called vertices) are, respectively,  
$$[0,3^{-j}], [3^{-j},2 \cdot 3^{-j}],\ldots,[1-3^{-j}, 1]\quad\mbox{and}\quad [0], [3^{-j}],\ldots,[1-3^{-j}], [1].$$
We denote by $I(m,j)$  the cell complex on $I^m$: 
$$I(m,j)=I(1,j)\otimes\ldots \otimes I(1,j)\quad (\mbox{$m$ times}).$$
Then $\alpha=\alpha_1 \otimes \cdots\otimes \alpha_m$ is a $q$-cell of $I(m,j)$ if and only if $\alpha_i$ is a cell
of $I(1,j)$ for each $i$, and $\sum_{i=1}^m {\rm dim}(\alpha_i) =q$. We often abuse notation by identifying  a $q$-cell $\alpha$ with its support: $\alpha_1 \times \cdots \times \alpha_m \subset I^m$. 

The cube complex $X(j)$ is the union of all cells of $I(m,j)$ whose support is contained in some cell of $X$. We use the notation $X(j)_q$ to denote the set of all $q$-cells in $X(j)$. Two vertices $x, y\in X(j)_0$ are {\em adjacent} if they belong to a common cell in $X(j)_1$.

Given $i,j\in \N$ we define ${\bf n}(i,j):X(i)_0\rightarrow X(j)_0$ so that ${\bf n}(i,j)(x)$ is the element in $X(j)_0$ that is closest to $x$ (see \cite[page 141]{pitts} or \cite[Section 7.1]{marques-neves} for a precise definition).  
 
 Given a map $\phi:X(j)_0\rightarrow  \mathcal{Z}_n(M;\Z_2)$, we define the {\em fineness} of $\phi$ to be
$${\bf f}(\phi)=\sup\left\{{\bf M}(\phi(x)-\phi(y)): x,y \mbox{ adjacent vertices in } X(j)_0\right\}.$$
The reader should think of the notion of fineness as being a discrete  measure of continuity with respect to the mass norm. 

\subsection{Homotopy notions}
Let $\phi_i:X(k_i)_0\rightarrow  \mathcal{Z}_n(M;\Z_2)$, $i=1,2$. We say that $\phi_1$ is {\it $X$-homotopic to $\phi_2$ in $\mathcal{Z}_n(M;{\bf M};\Z_2)$ with fineness $\delta$}  if we can find  $k\in \N$ and a map
$$\psi: I(1,k)_0\times X(k)_0\rightarrow  \mathcal{Z}_n(M;\Z_2)$$
such that 
\begin{itemize}
\item[(i)] ${\bf f}(\psi)<\delta;$
\item[(ii)] if $i=1,2$ and $x\in X(k)_0$, then
$$\psi([i-1],x)=\phi_i({\bf n}(k,k_i)(x)).$$
\end{itemize}

 Instead of considering continuous maps from $X$ into $\mathcal{Z}_n(M;{\bf M};\Z_2)$, the Almgren-Pitts theory deals with  sequences of discrete maps into  $\mathcal{Z}_n(M;\Z_2)$ with finenesses tending to zero. 
 
\subsection{Definition}\label{homotopy.sequence.phi}
 An $$\mbox{{\it $(X,{\bf M})$-homotopy sequence of mappings into $\mathcal{Z}_n(M;{\bf M};\mathbb{Z}_2)$}}$$ is a sequence of mappings $S=\{\phi_i\}_{i\in \N}$,
$$\phi_i:X(k_i)_0\rightarrow \mathcal{Z}_n(M;\mathbb{Z}_2),$$
such that $\phi_i$ is $X$-homotopic to $\phi_{i+1}$ in  $\mathcal{Z}_n(M;{\bf M};\mathbb{Z}_2)$ with fineness $\delta_i$ and
\begin{itemize}
\item[(i)] $\lim_{i\to\infty} \delta_i=0$;
\item[(ii)]$\sup\{{\bf M}(\phi_i(x)):x\in X(k_i)_0, i\in \N\}<+\infty.$
\end{itemize}

\medskip

The  next definition explains what  it means for two distinct homotopy sequences  of mappings into  $\mathcal{Z}_n(M;{\bf M};\Z_2)$ to  be homotopic.

\subsection{Definition}
Let $S^1=\{\phi^1_i\}_{i\in \N}$ and $S^2=\{\phi^2_i\}_{i\in \N}$  be $(X,{\bf M})$-homotopy sequences  of mappings into  $\mathcal{Z}_n(M;{\bf M};\mathbb{Z}_2)$. We say  that {\it $S^1$ is homotopic with $S^2$} if there exists a sequence $\{\delta_i\}_{i\in \N}$ such that
 \begin{itemize}
\item $\phi^1_i$  is $X$-homotopic to $\phi^2_i$ in  $\mathcal{Z}_n(M;{\bf M};\mathbb{Z}_2)$ with fineness $\delta_i$;
\item $\lim_{i\to\infty} \delta_i=0.$
 \end{itemize}
 
 \medskip
 
  The relation ``is homotopic with'' is an equivalence relation on the set of all $(X,{\bf M})$-homotopy sequences  of mappings into  $\mathcal{Z}_n(M;{\bf M};\mathbb{Z}_2)$. We call  the equivalence class 
of any such sequence  an {\it $(X,{\bf M})$-homotopy class of mappings into $\mathcal{Z}_n(M;{\bf M};\mathbb{Z}_2)$}. We denote 
by $[X,\mathcal{Z}_n(M;{\bf M};\mathbb{Z}_2)]^{\#}$  the set of all equivalence classes.

The definitions of homotopy  for sequences of  discrete maps whose  finenesses  are measured with respect to the flat metric, instead of the mass norm, are entirely analogous. These are discrete analogues 
of the usual notions of homotopy for continuous maps $\Phi:X\rightarrow \mathcal{Z}_n(M;\mathbb{Z}_2)$. 

\subsection{Width}
Given  $\Pi \in  [X,\mathcal{Z}_n(M;{\bf M};\mathbb{Z}_2)]^{\#}$, let
$$ {\bf L}: \Pi\rightarrow [0,+\infty]$$
be defined by
$${\bf L}(S)=\limsup_{i\to\infty}\max\{{\bf M}(\phi_i(x)):x\in \mathrm{dmn}(\phi_i)\},\quad\mbox{where }S=\{\phi_i\}_{i\in \N}.$$
Note that ${\bf L}(S)$ is the discrete replacement for  the maximum area of a continuous map into $\mathcal{Z}_n(M;{\bf M};\mathbb{Z}_2)$.

Given $S=\{\phi_i\}_{i\in \N}\in \Pi$, we also consider the compact subset ${\bf K}(S)$ of $\mathcal{V}_n(M)$ given by
\begin{multline*}
{\bf K}(S)=\{V:V=\lim_{j\to\infty}|\phi_{i_j}(x_j)|\mbox{ as varifolds, for some increasing}\\
\mbox{sequence }\{i_j\}_{j\in \N} \mbox{ and }x_j\in \mathrm{dmn}(\phi_{i_j})\}.
\end{multline*}
This is the discrete replacement for the image of a  continuous map into $\mathcal{Z}_n(M;{\bf M};\mathbb{Z}_2)$.

\subsection{Definition}The {\it width} of $\Pi$ is defined by
$${\bf L}(\Pi)=\inf\{{\bf L}(S):S\in \Pi\}.$$

We say  $S\in \Pi$ is a {\it critical sequence} for $\Pi$ if $${\bf L}(S)={\bf L}(\Pi).$$
The {\em critical set}  ${\bf C}(S)$ of a critical sequence   $S\in \Pi$ is given by 
$${\bf C}(S)={\bf K}(S)\cap\{V: ||V||(M)={\bf L}(S)\}.$$

Consider $\Pi \in [X,\mathcal{Z}_n(M;{\bf M};\mathbb{Z}_2)]^{\#}$.  {The next proposition states that tight critical sequences always exist}.

\subsection{Proposition}\label{pulltight}\textit{ Suppose  $\partial M=\emptyset$.
There exists  a critical sequence $S^* \in \Pi$. Moreover, for each  critical sequence $S^*\in \Pi$ there exists a critical sequence $S\in\Pi$ such that
\begin{itemize}
\item ${\bf C}(S)\subset {\bf C}(S^*)$;
\item every $\Sigma\in {\bf C}(S)$ is  a stationary varifold.
\end{itemize}
}
The sequence $S$ is obtained from a pull-tight procedure applied to  $S^*$. The proof is essentially the same of Theorem 4.3 of \cite{pitts} (see also Section 15 of \cite{marques-neves}).

\subsection{Almost minimizing varifold}
In order to explain  the  regularity part of the Almgren-Pitts min-max theory, we need to introduce the notion of an almost minimizing varifold. 

\subsection{Definition} A varifold $V\in \mathcal V_{n}(M)$ is {\em $\Z_2$ almost minimizing} in an open set $U\subset M$ if  for every $\varepsilon>0$ we can find $\delta>0$ and  
$$T\in \mathcal Z_{n}(M, M\setminus U;\Z_2),$$ with ${\bf F}_U(V,|T|)<\varepsilon$ and such that the following property holds true:

if $\{T_i\}_{i=0}^q$ is a sequence in $\mathcal Z_{n}(M, M\setminus U;\Z_2)$ with 
\begin{itemize}
\item $T_0=T$ and $\mbox{spt}(T-T_i)\subset U$ for all $i=1,\ldots, q;$
\item ${\bf M}(T_i-T_{i-1})\leq \delta$ for all $i=1,\ldots, q;$
\item ${\bf M}(T_i)\leq {\bf M}(T)+\delta$ for all $i=1,\ldots, q;$
\end{itemize}
 then ${\bf M}(T_q)\geq {\bf M}(T)-\varepsilon$.
\vskip 0.05in

Loosely speaking this is saying that every  deformation of  $V\in\mathcal V_{n}(M)$ that is supported in $U$ and that decreases the area by more than $\varepsilon$ must pass through a stage where the area is increased by more than  $\delta$.

\medskip

Given real numbers $0<s<r$, let $A(p,s,r)=\{ x\in \mathbb{R}^L: s<|x-p|< r\}$.

\subsection{Definition} A varifold $V \in \mathcal{V}_n(M)$ is {\it $\mathbb{Z}_2$ almost minimizing in annuli} if 
for each $p \in M$, there exists $r=r(p)>0$ such that $V$ is $\mathbb{Z}_2$ almost minimizing in $M \cap A (p,s,r)$ for all $0<s<r$. 

\medskip

If $V \in \mathcal{V}_n(M)$ is stationary in $M$ and $\mathbb{Z}_2$ almost minimizing in annuli, then $V \in \mathcal{IV}_n(M)$ by Theorem 3.13 of \cite{pitts}.

The regularity of almost minimizing integral varifolds was first done by Pitts in \cite[Section 7]{pitts} when $n\leq 5$, and then extended by Schoen and Simon to every dimension by allowing a singular set of codimension at least 7 \cite[Theorem 4]{schoen-simon}.  Schoen and Simon work with integer coefficients but, as we explain below, the arguments extend to $\Z_2$ coefficients also.
 
\subsection{Theorem}\label{regularity.thm} {\em  Suppose $n \leq 6$,  $\partial M=\emptyset$, and let $V \in \mathcal{IV}_{n}(M)$ be a nontrivial integral varifold that is both  stationary in $M$ and $\mathbb{Z}_2$ almost minimizing in annuli. Then $V$ is the varifold of a smooth, closed, embedded minimal hypersurface, with possible multiplicities.} 
\medskip
\begin{proof}
Let  $\mathcal{A}$ be the collection of all nontrivial $V \in \mathcal{IV}_{n}(M)$ that are stationary in $M$ and $\mathbb{Z}_2$ almost minimizing in annuli.

It follows from the work of Pitts  in \cite[Theorem 3.11]{pitts} that for any  $p \in {\rm spt}||V||$, we can find $r(p)>0$ such that for any $0<s<t<r(p)$ there
exists a replacement varifold $V^* \in \mathcal{A}$ with the properties:
\begin{itemize}
\item[(i)]  $||V^*||(M) = ||V||(M)$,
\item[(ii)] $V^* \llcorner G_{n}(M \setminus \overline{A}(p,s,t))=V\llcorner G_{n}(M \setminus \overline{A}(p,s,t))$,
\item[(iii)] $V^* \llcorner G_{n}(M \cap A(p,s,t))=(\lim_{j \to \infty} |T_j|) \llcorner G_{n}(M \cap A(p,s,t))$,
\end{itemize}
with $\{ T_j\} \subset {\bf I}_{n}(M, \mathbb{Z}_2)$, $\{ {\bf M}(T_j)\}$ bounded independently of $j$,
${\rm spt}(\partial T_j) \cap A(p,s,t) = \emptyset$, $T_j$ locally area minimizing in $M \cap A(p,s,t)$ and $|T_j|$
stable in $M \cap A(p,s,t)$. By choosing $r(p)$ sufficiently small, we also get that $M \cap A(p,s,t)$ is simply connected
for every $0<s<t<r(p)$.

 It follows from the regularity theory for  area minimizing mod 2 flat chains in \cite[Regularity Theorem 2.4]{morgan} (all conditions are satisfied by Remark 1 in \cite[page 249]{morgan}) 
 that there exists a smooth minimal hypersurface $\Sigma_j$ properly embedded  in $A(p,s,t)$ such that
$$({\rm spt \ }T_j) \cap A(p,s,t) = \overline{\Sigma}_j \cap A(p,s,t).$$
Since $M \cap A(p,s,t)$ is simply connected, we have that $\Sigma_j$ is orientable for each $j$. Therefore 
$${\rm spt \ }||V^*|| \cap A(p,s,t)  = \overline{\Sigma} \cap A(p,s,t) ,$$ where $\Sigma$ is an orientable stable smooth minimal hypersurface
exactly like in Schoen-Simon \cite[page 789]{schoen-simon}. From this point on, the proof that ${\rm spt}||V||$ is a smooth embedded minimal hypersurface 
proceeds just like in the proof of \cite[Theorem 4]{schoen-simon}.

\end{proof}

\subsection{Existence of almost minimizing varifolds}\label{am.varifolds}
The existence of almost minimizing varifolds is achieved in Theorem 4.10 of Pitts book \cite{pitts} through a combinatorial argument. This was inspired by a previous construction of Almgren \cite{almgren-varifolds} and is a crucial part of the Almgren--Pitts theory. The idea is that if $S$ is a homotopy sequence of maps such that every element in ${\bf C}(S)$ is stationary and no element in ${\bf C}(S)$ is almost minimizing in annuli, then the combinatorial arguments in \cite[page $165$--page $174$]{pitts} give a new homotopy sequence $S^*$ homotopic with $S$ such that ${\bf L}(S^*)<{\bf L}(S)$.

For the application we have in mind, the discrete maps in our sequence  are not defined on the whole grid $I(m,k_i)_0$ but only on the vertices of a subcomplex $Y_i$ of $I(m,k_i)$. Nonetheless, Pitts arguments immediately adapt to this setting and give the result that we now state in a precise way.

Consider a sequence of cubical subcomplexes $Y_i$ of $I(m,k_i)$, with $k_i \rightarrow \infty$,  and a sequence $S=\{\varphi_i\}$  of maps 
$$\varphi_i:(Y_i)_0 \rightarrow \mathcal{Z}_n(M;\Z_2),$$
with finesses $\delta_i$ tending to zero. Similarly as before, we define
$${\bf L}(S)=\limsup_{i\to\infty}\max\{{\bf M}(\varphi_i(x)):x\in \mathrm{dmn}(\varphi_i)\},$$
\begin{multline*}
{\bf K}(S)=\{V \in \mathcal{V}_n(M):V=\lim_{j\to\infty}|\varphi_{i_j}(x_j)|\mbox{ as varifolds, for some increasing}\\
\mbox{sequence }\{i_j\}_{j\in \N} \mbox{ and }x_j\in \mathrm{dmn}(\varphi_{i_j})\}.
\end{multline*}
and
$${\bf C}(S)={\bf K}(S)\cap\{V: ||V||(M)={\bf L}(S)\}.$$

If $Y$ is a subcomplex of $I(m,k)$, then similarly as before we define the cube subcomplex $Y(l)$ to be the the union of all cells of $I(m,k+l)$ whose support is contained in some cell of $Y$. The same notion of homotopy with fineness $\delta$ applies to maps $\phi_1:Y(l_1) \rightarrow \mathcal{Z}_n(M;\Z_2)$ and $\phi_2:Y(l_2) \rightarrow \mathcal{Z}_n(M;\Z_2)$.

\subsection{Theorem} \label{combinatorial.pitts}{\em Suppose  $\partial M=\emptyset$. Let $S=\{\varphi_i\}$ be as above, and such that every $V \in {\bf C}(S)$ is stationary in $M$.  If no element $V \in {\bf C}(S)$ is $\mathbb{Z}_2$ almost minimizing in annuli, then there exists a sequence $S^*=\{\varphi_i^*\}$ of maps 
$$\varphi_i^*: Y_i(l_i)_0 \rightarrow \mathcal{Z}_n(M;\Z_2),$$
for some $l_i \in \mathbb{N}$, such that:
\begin{itemize}
\item $\varphi_i$ and $\varphi_i^*$ are homotopic to each other with finesses that tend to zero as $i \rightarrow \infty$,
\item ${\bf L}(S^*)=\limsup_{i\to\infty}\max\{{\bf M}(\varphi_i^*(y)):y\in Y_i(l_i)_0\}<{\bf L}(S).$
\end{itemize} 
}
\medskip
Given $\Pi \in [X,\mathcal{Z}_n(M;{\bf M};\mathbb{Z}_2)]^{\#}$ we can apply this result to the critical sequence given by Proposition \ref{pulltight} and obtain the following simple extension of Theorem 4.10 in \cite{pitts}.
\subsection{Theorem}\label{pitts.min.max} {\em   
Suppose  $\partial M=\emptyset$, and let $\Pi \in [X,\mathcal{Z}_n(M;{\bf M};\mathbb{Z}_2)]^{\#}$. Then there exists an integral varifold $V \in \mathcal{IV}_n(M)$ 
such that the following three statements are true:
\begin{itemize}
\item [(1)] $||V||(\mathbb{R}^L) = {\bf L}(\Pi)$,
\item[(2)] $V$ is stationary in $M$,
\item [(3)] $V$ is $\mathbb{Z}_2$ almost minimizing in annuli.
\end{itemize}
Moreover, if $S^*$ is a critical  sequence for $\Pi$   then we can choose $V \in {\bf C}(S^*)$.
}

\section{Almgren's isomorphism and Interpolation results}\label{almgren.homo}

We describe some of the maps defined by Almgren in \cite[Section 3]{almgren}. There he uses  
integer coefficients and the unit interval $[0,1]$ as the parameter space, but everything extends to the setting
of  $\Z_2$ coefficients and of maps parametrized by the circle  $S^1$ instead.

Almgren associates to every continuous map in the flat topology $\Phi$ from $S^1$ into $\mathcal{Z}_n(M;\Z_2)$ (or $\mathcal{Z}_n(M,\partial M;\Z_2)$), an element $F(\Phi)$ in $H_{n+1}(M,\Z_2)$ (or $H_{n+1}(M,\partial M;\Z_2))$ such that $F(\Phi)=0$ if and only if $\Phi$ is homotopically trivial. He also provides equivalent constructions for discrete maps. We need both aspects of the theory and so we review his constructions and the interpolation results needed to make sure that one can move consistently from continuous maps to discrete maps.

\subsection{Discrete setting}\label{almgren.homo.discrete}
Suppose  we have a  map $$\phi:I(1,k)_0\rightarrow \mathcal{Z}_n(M,\partial M;\Z_2),$$ with  $\phi([0])=\phi([1])$ and so that
\begin{equation*}
\mathcal{F}(\phi(a_j),\phi(a_{j+1}))\leq \nu_{M,\partial M}\quad\mbox{for all }j=0,\dots, 3^k-1,
\end{equation*}
where $a_j=[j 3^{-k}]$ and $\nu_{M,\partial M}$,  defined in \cite[Theorem 2.4]{almgren}, is a small positive constant that depends only on $M$. This condition ensures the existence  of a constant $\rho=\rho(M)\geq 1$ and of {\em isoperimetric choices} $A_j\in {\bf I}_{n+1}(M;\Z_2)$ such that
$$\partial A_j-(\phi (a_{j+1})-\phi(a_j))\in {\bf I}_n(\partial M;\Z_2)\mbox{ and }{\bf M}(A_j)< \rho\mathcal{F}(\phi(a_j),\phi(a_{j+1}))$$
for all $j=0,\dots, 3^k-1.$
Hence $\sum_{j=0}^{3^k-1}A_j\in \mathcal{Z}_{n+1}(M,\partial M;\Z_2)$ and therefore it  defines a relative homology class (see \cite[Section 4.4]{federer}):
$$ F_{M,\partial M}^{\#}(\phi)=\left[\sum_{j=0}^{3^k-1}A_j\right]\in H_{n+1}(M,\partial M;\Z_2).$$ 
The following simple lemma shows that the  isoperimetric choice is unique.
\subsection{Lemma}\label{constancy.lemma}{\em The constant $\nu_{M,\partial M}$ can be chosen so that if $C_j\in{\bf I}_{n+1}(M;\Z_2)$ has
$${\bf M}(C_j)\leq \nu_{M,\partial M}\quad\mbox{and}\quad \partial C_j-(\phi (a_{j+1})-\phi(a_j))\in {\bf I}_{n}(\partial M;\Z_2),$$
then $A_j=C_j$.}
\begin{proof}
We have ${\rm spt}(\partial (A_j-C_j))\subset \partial M$ and so, by the Constancy Theorem \cite[Theorem 26.27]{simon}, we have $A_j-C_j=kM$ for some $k\in\{0,1\}$. 
Furthermore
$${\bf M}(A_j)\leq \rho  \mathcal{F}(\phi(a_j),\phi(a_{j+1})) \leq \rho\nu_{M,\partial M}.$$
Thus ${\bf M}(A_j-C_j)\leq (\rho+1)\nu_{M,\partial M}$.  The result follows if  $(\rho+1)\nu_{M,\partial M}$ is strictly smaller than ${\bf M}(M)$. 
\end{proof}
The work of Almgren \cite{almgren} shows that if another map
$$\phi':I(1,k')_0\rightarrow \mathcal{Z}_n(M,\partial M;\Z_2),$$ with  $\phi'([0])=\phi'([1])$, is homotopic to $\phi$ in the discrete sense, with fixed boundary values, and with 
fineness in the flat topology smaller than $\nu_{M,\partial M}$, then 
\begin{equation}\label{almgren.homo.boundary}
F_{M,\partial M}^{\#}(\phi)=F_{M,\partial M}^{\#}(\phi').
\end{equation}

\subsection{Continuous setting} Assume $\partial M=\emptyset$ for simplicity.
Given a continuous map in the flat topology $$\Phi:S^1\rightarrow \mathcal{Z}_n(M;\Z_2),$$  
 we can take $k$ sufficiently large so that, 
\begin{equation}\label{almgren.constant}
\mathcal{F}(\Phi(e^{2\pi ix}),\Phi(e^{2\pi iy}))\leq \nu_{M}\quad\mbox{for all }x, y \mbox{ in a common cell of  }I(1,k).
\end{equation}

If $\phi:I(1,k)_0\rightarrow \mathcal{Z}_n(M;\Z_2)$ is given by $\phi([x])=\Phi(e^{2\pi ix})$, we can define
$$F_{M}(\Phi)=  F_{M}^{\#}(\phi)\in H_{n+1}(M,\Z_2).$$
We have that the homology class $F_M(\Phi)$ does not depend on $k$,  provided condition \eqref{almgren.constant} is satisfied, and that  $$F_{M}(\Phi)=F_{M}(\Phi')$$
for any continuous map $\Phi':S^1\rightarrow \mathcal{Z}_n(M;\Z_2)$ in the homotopy class of $\Phi$.  {Moreover, Almgren's work  \cite{almgren} also shows  that the  induced map
$$F_M:\pi_1(\mathcal{Z}_n(M;\Z_2))\rightarrow H_{n+1}(M;\Z_2),\quad[\Phi]\mapsto [F_M(\Phi)] $$
is an isomorphism.}

\subsection{Definition} A continuous map  in the flat topology $\Phi:S^1\rightarrow \mathcal{Z}_n(M;\Z_2)$ with $F_{M}(\Phi)\neq 0$ is called a  {\em  sweepout } of $M$. If $F_{M}(\Phi)= 0$, we say $\Phi$ is {\em trivial}.

\medskip
The next proposition follows from the work of Almgren \cite{almgren} and its proof is left to Appendix \ref{A.section}.

\subsection{Proposition}\label{close.implies.homotopic} \textit{Let $Y$ be a cubical subcomplex of some $I(m,l)$.  There exists $\delta=\delta(M,m)>0$ with the following property:}

\textit{ 
If $\Phi_1,\Phi_2:Y \rightarrow \mathcal{Z}_n(M;\Z_2)$ are continuous maps in the flat topology such that 
$$\sup\{\mathcal{F}(\Phi_1(y),\Phi_2(y)):y\in Y\} < \delta,$$ then $\Phi_1$ is homotopic to $\Phi_2$ in the flat topology. }
\medskip 

One immediate consequence is the following corollary:

\subsection{Corollary}\label{locally.trivial} {\em Let $\mathcal{T}$ be a finite subset of $\mathcal{Z}_n(M;\Z_2)$. If $\varepsilon>0$
is sufficiently small, depending on $\mathcal{T}$, then every map $\Phi:S^1\rightarrow \mathcal{Z}_n(M;\Z_2)$ with
$$\Phi(S^1) \subset B_\varepsilon^\mathcal{F}(\mathcal{T}) = \{T \in \mathcal{Z}_n(M;\Z_2): \mathcal{F}(T,\mathcal{T})< \varepsilon\}$$ is trivial. 
} 

\begin{proof}
Let $d=\min \{\mathcal{F}(S,T): S,T\in \mathcal{T}, S\neq T \}$ and set $\varepsilon=\min\{\delta,d/3\}$, where $\delta$ is given by Proposition \ref{close.implies.homotopic}. 

The fact that $\Phi(S^1) \subset B_\varepsilon^\mathcal{F}(\mathcal{T})$ implies that $\Phi(S^1) \subset B_\varepsilon^\mathcal{F}(T)$ for some $T\in \mathcal T$. Thus, Proposition \ref{close.implies.homotopic} implies that $\Phi$ is homotopic to a constant map $\Phi'$ and so $F_M(\Phi)=F_M(\Phi')=0$.




\end{proof}
\subsection{Interpolation results}
Given a continuous map $\Phi:X\rightarrow {\mathcal Z}_n(M;\Z_2)$, with respect to the flat topology, we say that $\Phi$  {\em has no concentration of mass} if
$$\lim_{r\to 0} \sup\{||\Phi(x)||(B_r(p)):x\in X, p\in M\}=0.$$
This is a mild technical condition  which is satisfied by all maps we construct in this paper.
\subsection{Lemma}\label{mass.concentration.lemm}{\em If  $\Phi:X\rightarrow {\mathcal Z}_n(M;{\bf M};\Z_2)$ is continuous  in the mass norm, then 
$$\sup\{{\bf M}(\Phi(x)):x\in X\}<+\infty$$
and $\Phi$ has no concentration of mass.}
\begin{proof}

Choose $\delta >0$. Given $p\in M$ and $x \in X$, there is $r=r(p,x)>0$ and $U_{(p,x)} \subset X$ an open neighborhood of $x$ so that
$$
||\Phi(y)||(B_r(p))< \delta 
$$
for all $y \in U_{(p,x)}$.

By compactness, we can select a finite covering $\{B_{r_k}(p_k) \times U_{(p_k,x_k)}\}_{k=1}^N$ of $M \times X$, where $r_k=r(p_k,x_k)/2$. If $\overline{r}=\min\{r_k\}_{k=1}^N$, then 
$$
||\Phi(x)||(B_{\overline{r}}(p))< \delta 
$$
for all $(p,x) \in M \times X$ and the result follows.


\end{proof}

The next  theorem follows from Theorem 13.1 in \cite{marques-neves} and its purpose is to construct a $(X,{\bf M})$-homotopy sequence of mappings out of a continuous map in the flat topology with no concentration of mass.
 
\subsection{Theorem}\label{continuous.discrete} {\em Let $\Phi:X\rightarrow \mathcal{Z}_n(M;\Z_2)$ be a continuous map in the flat topology  that has no concentration of  mass. 
There exist a sequence of maps
$$\phi_i:X(k_i)_0 \rightarrow \mathcal{Z}_n(M;\Z_2),$$
with $k_i<k_{i+1}$, and a
sequence of positive numbers $\{\delta_i\}_{i\in\N}$ converging to zero such that
\begin{itemize}
\item[(i)] $$S=\{\phi_i\}_{i\in\N}$$ is an $(X,{\bf M})$-homotopy sequence of mappings into $\mathcal{Z}_n(M;{\bf M};\Z_2)$ with ${\bf f}(\phi_i)<\delta_i$;
\item[(ii)] $$\sup\{\mathcal F(\phi_i(x)-\Phi(x)): x\in X(k_i)_0\}\leq \delta_i;$$
\item[(iii)]$$\sup\{{\bf M}(\phi_i(x)): x\in X(k_i)_0\}\leq \sup\{{\bf M}(\Phi(x)): x\in X\}+\delta_i.$$
\end{itemize}
}

The next theorem follows from Theorem 14.1 in \cite{marques-neves} and its purpose is to construct a continuous map in the mass norm out of a discrete map with small fineness.
\subsection{Theorem}\label{discrete.continuous} {\em There exist positive constants $C_0=C_0(M,m)$ and $\delta_0=\delta_0(M)$  so that if $Y$ is a cubical subcomplex of $I(m,k)$ and 
$$\phi:Y_0\rightarrow \mathcal{Z}_n(M;\Z_2)$$
has ${\bf f}(\phi)<\delta_0$, then there exists   a map 
$$ \Phi:Y\rightarrow \mathcal{Z}_n(M;{\bf M};\Z_2)$$
continuous in the mass norm and satisfying  
\begin{itemize}
\item[(i)]$\Phi(x)=\phi(x)$ for all $x\in Y_0$;
\item[(ii)] if $\alpha$ is some $j$-cell in $Y_j$, then $\Phi$ restricted to $\alpha$ depends only on the values of $\phi$ assumed on the vertices of  $\alpha$;
\item[(iii)]
$$\sup\{{\bf M}(\Phi(x)-\Phi(y)): x,y\mbox{ lie in a common cell of } Y\}\leq C_0{\bf f}(\phi).$$
\end{itemize}
}
We call the map $\Phi$ given by Theorem \ref{discrete.continuous} \textit{the Almgren extension} of $\phi$. The next proposition shows that the Almgren extension preserves the homotopy classes.

\subsection{Proposition}\label{discrete.homotopy.continuous}{\em Let $Y$ be a cubical subcomplex of $I(m,k)$. There exists $\eta=\eta(M,m)>0$ with the following property:

 If $\phi_1:Y(l_1)_0 \rightarrow \mathcal{Z}_n(M;\Z_2)$ is homotopic to $\phi_2:Y(l_2)_0 \rightarrow \mathcal{Z}_n(M;\Z_2)$  with fineness smaller than $\eta$, then the Almgren extensions
 $$\Phi_1, \Phi_2:Y\rightarrow \mathcal{Z}_n(M;{\bf M};\Z_2)$$  of $\phi_1, \phi_2$, respectively, are homotopic to each other in the flat topology.
}
\begin{proof}Set $\eta=\delta/(2C_0)$, where $\delta$ and $C_0$ are given by Proposition \ref{close.implies.homotopic} and Theorem \ref{discrete.continuous}, respectively.

By assumption, we can find  $l\in \N$ and a map
$$\psi: I(1,k+l)_0\times Y(l)_0\rightarrow  \mathcal{Z}_n(M;\Z_2)$$
with ${\bf f}(\psi)<\eta$ and such that if $i=1,2$ and $y\in Y(l)_0$, then
$$\psi([i-1],y)=\phi_i({\bf n}(k+l,k+l_i)(y)).$$

For $i=1,2$, let $\phi_i':Y(l)_0 \rightarrow \mathcal{Z}_n(M;\Z_2)$ be given by $\phi_i'(y)=\psi([i-1],y)$ and let 
$\Phi_i':Y\rightarrow \mathcal{Z}_n(M;{\bf M};\Z_2)$ be the Almgren extension of $\phi_i'$.

 By Theorem \ref{discrete.continuous}, it follows that ${\bf M}(\Phi_i(y),\Phi_i'(y))\leq 2C_0\eta\leq \delta$ for every $y \in Y$ and so Proposition \ref{close.implies.homotopic} implies that $\Phi_i$ is homotopic to $\Phi_i'$ in the flat topology, for each $i=1,2$. 
The Almgren extension of $\psi$ to $I\times Y$ is a homotopy between $\Phi_1'$ and $\Phi_2'$ and this implies the result.

\end{proof}
We end this section with the following corollary.

\subsection{Corollary}\label{discrete.continuous.homotopy}{\em Let $S=\{\phi_i\}_{i\in \N}$ and $S'=\{\phi'_i\}_{i\in \N}$  be $(X,{\bf M})$-homotopy sequences  of mappings into  $\mathcal{Z}_n(M;{\bf M};\mathbb{Z}_2)$ such that  $S$ is homotopic with $S'$.
\begin{itemize}
 \item[(i)] The Almgren extensions of $\phi_i$, $\phi_i'$:
 $$\Phi_i, \Phi_i':X\rightarrow \mathcal{Z}_n(M;{\bf M};\Z_2),$$ respectively, are homotopic to each other in the flat topology for  sufficiently large $i$.
 \item[(ii)] {If $S$ is given by Theorem \ref{continuous.discrete} (i) applied to $\Phi$, where  $\Phi:X\rightarrow \mathcal{Z}_n(M;\Z_2)$ is a continuous map in the flat topology with no concentration of  mass,} then $\Phi_i$ is homotopic to $\Phi$ in the flat topology  for every sufficiently large $i$.
 Moreover,
$$
\limsup_{i\rightarrow \infty} \sup\{ {\bf M}(\Phi_i(x)):x\in X\}={\bf L}(S) \leq  \sup\{ {\bf M}(\Phi(x)):x\in X\}.
$$
 
\end{itemize}
}
\begin{proof}
Property (i) follows immediately from Proposition \ref{discrete.homotopy.continuous} and the definition of homotopy between
sequences of mappings into $\mathcal{Z}_n(M;{\bf M};\mathbb{Z}_2)$.

From Theorem \ref{continuous.discrete} (i) and (ii), and Theorem \ref{discrete.continuous} (i) and (iii)
$$\lim_{i\to\infty}\sup\{\mathcal{F}(\Phi_i(x),\Phi(x)):x\in X\}=0$$
and thus, by  Proposition \ref{close.implies.homotopic}, $\Phi_i$ is homotopic to $\Phi$ in the flat topology  for all $i$ sufficiently large.
The statement about the supremum of the masses follows from Theorem \ref{continuous.discrete} (i) and (iii), and Theorem \ref{discrete.continuous} (i) and (iii).
\end{proof}


\section{Min-max Families}\label{minmax.families}

In this section we denote by $X$ a  cubical subcomplex of  $I^m=[0,1]^m$, for some $m$. 

The Almgren isomorphism $F_M$ establishes an isomorphism between $\pi_1(\mathcal{Z}_n(M;\Z_2))$ and $H_{n+1}(M;\Z_2) {=} \Z_2$. Hence $$H^1( \mathcal{Z}_n(M;\Z_2); \Z_2) {=} \Z_2$$ with a generator $\overline{\lambda}$. Denote by   $\bar \lambda^p$  the cup product of $\bar \lambda$ with itself $p$ times.

\subsection{Definition}\label{defi.pclass}{A continuous map $\Phi:X\rightarrow \mathcal{Z}_n(M;\Z_2)$ is a {\em  $p$-sweepout} if
$$
\Phi^*(\bar \lambda^p) \neq 0 \in H^p(X;\Z_2).
$$
This is equivalent  to say that there exists $\lambda \in H^1(X;\Z_2)$ such that:
\begin{itemize}
\item[(i)] for any cycle  $\gamma:S^1 \rightarrow X$,  we have $\lambda(\gamma)\neq 0$ if and only if  $\Phi\circ\gamma:S^1 \rightarrow \mathcal{Z}_n(M;\Z_2)$ is a  sweepout;
\item[(ii)] the cup product $\lambda^p=\lambda\smile\ldots\smile\lambda$ is nonzero in $H^p(X;\Z_2)$. 
\end{itemize}}
\subsection{Remark:}  \begin{enumerate}
\item A continuous map in the flat topology that is   homotopic to a $p$-sweepout is also a $p$-sweepout.
\item If $\gamma$, $\gamma'$ are homotopic to each other in $X$, then  $\Phi\circ\gamma$ is a  sweepout if and only if $\Phi\circ\gamma'$ is a sweepout. This will be useful to check condition (i) above in specific examples.
\end{enumerate}

\medskip

We say $X$ is {\em $p$-admissible} if there exists  a $p$-sweepout $\Phi:X\rightarrow \mathcal{Z}_n(M;\Z_2)$  that has no concentration of mass. The set of all $p$-sweepouts $\Phi$ that have no concentration of mass is denoted by $\mathcal P_p$. Note that two maps in  $\mathcal P_p$ can have different domains.

Similarly to Guth \cite[Appendix 3]{guth}, we define
\subsection{Definition}\label{p.width} The {\em $p$-width of $M$}  is $$\omega_p(M)=\inf_{\Phi \in \mathcal P_p}\sup\{{\bf M}(\Phi(x)): x\in {\rm dmn}(\Phi)\},$$
where ${\rm dmn}(\Phi)$ is the domain of $\Phi$.

\medskip
Notice that if a map $\Phi:X\rightarrow \mathcal{Z}_n(M;\Z_2)$ is a $p$-sweepout, then it also a $q$-sweepout
for every $q <p$. Hence $\omega_p(M) \leq \omega_{p+1}(M)$ for every $p \in \mathbb{N}$.

\subsection{Definition} Let $\Pi \in  [X,\mathcal{Z}_n(M;{\bf M};\mathbb{Z}_2)]^{\#}$. We say that $\Pi$ is a {\em{class of (discrete)  $p$-sweepouts}} if for any $S=\{\phi_i\}\in \Pi$, the Almgren extension $\Phi_i:X\rightarrow \mathcal{Z}_n(M;{\bf M};\Z_2)$ of $\phi_i$ is a $p$-sweepout for every sufficiently large $i$.

\subsection{Remark}\label{remark.detects} By Corollary \ref{discrete.continuous.homotopy} (i), it is enough to check that this is true for some $S=\{\phi_i\}\in \Pi$.

\medskip
The next lemma assures us that the discrete and  continuous definitions of a $p$-sweepout  are consistent.

\subsection{Lemma}\label{almgren.detects} \textit{Let \begin{itemize}
\item $\Phi:X\rightarrow \mathcal{Z}_n(M;\Z_2)$  be  a continuous map in the flat topology with no concentration of mass;
\item $S=\{\phi_i\}$ be the sequence of discretizations associated to $\Phi$ given by Theorem \ref{continuous.discrete} (i);
\item$\Pi$ be  the  $(X,{\bf M})$-homotopy class of mappings into $\mathcal{Z}_n(M;{\bf M};\mathbb{Z}_2)$ associated with $S=\{\phi_i\}$.
\end{itemize}
Then $\Phi \in \mathcal P_p$ is a $p$-sweepout if and only  if $\Pi$ is a {class of  $p$-sweepouts.}}
\begin{proof}
Denote by $\Phi_i$ the Almgren extension of $\phi_i$. The map $\Phi_i$ is continuous in the mass norm and hence it has no concentration of mass (Lemma \ref{mass.concentration.lemm}).
Since $\Phi_i$ is homotopic to $\Phi$ in the flat topology for all large $i$, by Corollary \ref{discrete.continuous.homotopy} (ii), the lemma follows at once.
\end{proof}

The same consistency {between discrete and continuous definitions}  also holds for the $p$-width.

\subsection{Lemma}\label{pwidth.discrete}\textit{Let $\mathcal D_p$ be the set of all {classes of $p$-sweepouts} $$\Pi \in  [X,\mathcal{Z}_n(M;{\bf M};\mathbb{Z}_2)]^{\#},$$ where $X$ is any $p$-admissible cubical subcomplex. Then
$$\omega_p(M)=\inf_{\Pi\in \mathcal D_p}{\bf L}(\Pi).$$}
\begin{proof}
We claim that for any  $p$-admissible $X$ and any {class of $p$-sweepouts} $\Pi\in[X,\mathcal{Z}_n(M;{\bf M};\mathbb{Z}_2)]^{\#}$, we have $\omega_p(M)\leq {\bf L}(\Pi)$.

Indeed, choose $S=\{\phi_i\}\in \Pi$ with $ {\bf L}(S)\leq {\bf L}(\Pi)+\varepsilon$ (with $\varepsilon>0$ arbitrary), and let $\Phi_i$ denote the Almgren extension of each $\phi_i$. We have by Theorem \ref{discrete.continuous} (i) and (iii) that
$$\omega_p(M)\leq \limsup_{i\to\infty}\sup\{{\bf M}(\Phi_i(x)):x\in X\}={\bf L}(S)\leq {\bf L}(\Pi)+\varepsilon.$$
By letting $\varepsilon$ tend to zero we obtain the desired claim.

Now, let  $\varepsilon>0$ and choose $\Phi\in \mathcal P_p$ with $$\sup\{{\bf M}(\Phi(x)):x\in {\rm dmn}(\Phi)\}\leq \omega_p(M)+\varepsilon.$$ Consider $S$ and $\Pi$ as in the statement of Lemma \ref{almgren.detects}. Then $\Pi$ is a {class of $p$-sweepouts} and from Theorem \ref{continuous.discrete}  (iii) we have
$${\bf L}(\Pi)\leq{\bf L}(S)\leq  \sup\{{\bf M}(\Phi(x)):x\in {\rm dmn}(\Phi)\}\leq \omega_p(M)+\varepsilon.$$
By letting $\varepsilon$ tend to zero and using the previous claim  we prove the lemma.
\end{proof}

{It is not clear a  priori whether the number $\omega_p(M)$ is equal to the width ${\bf L}(\Pi)$ of some class of $p$-sweepouts $\Pi.$ The next proposition analyzes the case where this is not true.}

\subsection{Proposition}\label{minmax.cohomology}\textit{Assume  $2\leq n\leq 6$. If there exists $p\in\N$ such that for all  $p$-admissible $X$ we have
$$\omega_p(M)<{\bf L}(\Pi)\quad\mbox{for  every { class of $p$-sweepouts} }\Pi \in  [X,\mathcal{Z}_n(M;{\bf M};\mathbb{Z}_2)]^{\#},$$
then there exist infinitely many distinct smooth closed minimal embedded hypersurfaces with  uniformly bounded area.}

\begin{proof} From Lemma \ref{pwidth.discrete} we can find sequences of $p$-admissible cubical subcomplexes $X_k$ and  of {classes of $p$-sweepouts} $\Pi_k\in  [X_k,\mathcal{Z}_n(M;{\bf M};\mathbb{Z}_2)]^{\#}$ such that
$${\bf L}(\Pi_1)> \dots > {\bf L}(\Pi_k) > {\bf L}(\Pi_{k+1})> \dots$$ and $$\lim_{k \rightarrow \infty} {\bf L}(\Pi_k)= \omega_p(M).$$
The combination of Theorem \ref{pitts.min.max} and Theorem \ref{regularity.thm}  implies ${\bf L}(\Pi_k) = ||V_k||(M)$ for some smooth closed embedded minimal hypersurface $V_k$, possibly disconnected and with integer multiplicities. The proposition follows.
\end{proof}

\section{Upper bounds}\label{guth}

The asymptotic behavior of the min-max volumes $\omega_p(M)$ as $p \rightarrow \infty$ has been studied previously by Gromov and Guth. In \cite{guth}, Guth uses a bend--and--cancel argument to prove the following result, {which was also proven by Gromov in \cite[Section 4.2.B]{gromov0}}.

\subsection{Theorem}\label{upper.bound}{\em   For each $p\in \N$, there exists a map  $$\Phi:\RP^p\rightarrow \mathcal{Z}_n(M;\Z_2)$$ that is continuous in the flat topology, has no concentration of mass and which is a $p$-sweepout ($\Phi \in \mathcal{P}_p$). Moreover, there exists a constant $C=C(M)>0$ so that $$\omega_p(M)\leq  \sup_{x\in \RP^p}{\bf M}(\Phi(x)) \leq C p^{\frac{1}{n+1}}$$
for every $p\in \mathbb{N}$.
}

\medskip

Guth proved this theorem  in \cite[Section 5]{guth} when the ambient space is a unit ball, but the arguments carry over  to the case when the ambient space is a closed manifold $M$. We present them here for convenience of the reader.

Any compact differentiable manifold can be triangulated. Therefore, by \cite[Chapter 4]{panov}, we can find an $(n+1)$-dimensional cubical subcomplex  $K$ of $I^m$ for some $m$, and a Lipschitz homeomorphism $G:K\rightarrow M$ such that  $G^{-1}: M \rightarrow K$ is also   Lipschitz. For each $k\in \N$, we denote by $c(k)\subset M$ the  image under $G$ of the set consisting of the centers of the cubes $\sigma\in K(k)_{n+1}$ (recall the definition of $K(k)_p$ in Section \ref{gmt}). In what follows we abuse notation and identify  cells in the subdivision $K(k)$ with their support.

We need to establish some preliminary results.  The first lemma follows  from the local description of a Morse function in terms of linear or quadratic functions and we leave its proof to the reader.

\subsection{Lemma}\label{morse.flat}{\em Let $f:M \rightarrow \mathbb{R}$ be a Morse function. Then the following properties are true:
\begin{itemize}
\item [(i)] the level set $\Sigma_t=\{x \in M:f(x)=t\}$ has finite $n$-dimensional Hausdorff measure for every $t \in \mathbb{R}$;
\item[(ii)] for every $\varepsilon>0$ and $x\in M$, there exists a radius $r>0$ such that
$$
\mathcal{H}^n(\Sigma_t \cap B_r(x)) < \varepsilon
$$
for all $t \in \mathbb{R}$;
\item[(iii)] for every $\varepsilon>0$ there exists $\delta>0$ such that 
$$
|b-a|<\delta \Longrightarrow \mathrm{vol}\left(f^{-1}([a,b])\right)< \varepsilon.
$$
\end{itemize}
}

The next lemma uses the embedding of $M$ into some $\R^L$ to produce a suitable Morse function.

\subsection{Lemma}\label{morse}{\em Fix $k \in \mathbb{N}$. For almost all  $v\in S^{L-1}=\{x\in \R^L:|x|=1\}$, we have that
\begin{itemize}
 \item[(i)] the function $f:M \rightarrow \mathbb{R}$, with $f(x)=\langle x, v\rangle$,  is Morse;
\item[(ii)] $f^{-1}(t)\cap c(k)$ contains at most one point for all $t\in \R$;
 \item[(iii)] no critical point of $f$ belongs to $c(k)$.
\end{itemize}
 }
\begin{proof}
By Sard's theorem, the  function $f_v(x)=\langle x, v\rangle$, $x \in M$, is Morse for all $v$ in an open  subset $A$ of $S^{N-1}$ with full measure.
Consider 
$$B=\{v\in S^{L-1}: \langle v,u-w\rangle \neq 0\mbox{ for all }u, w\in c(k)\mbox{ with }u\neq w\}.$$
Hence $B$ is an open set with full measure.  Given $x\in M$, let $T_x^\perp M$ be the orthogonal complement of $T_xM$ in $\mathbb{R}^L$. Then the set 
$$
C=\{v\in S^{L-1}: v \notin T_u^\perp M \mbox{ for all }u\in c(k)\}
$$
is also open with full measure. The properties (i), (ii) and (iii) are satisfied for every $v\in A \cap B \cap C$, an open set with full measure.
\end{proof}

Finally, to apply Guth's bend--and--cancel argument, we need a Lipschitz map homotopic to the identity that  maps  the complement of  a small neighborhood of $c(k)$ in $M$ into the $n$-skeleton $G(K(k)_n)$.

\subsection{Proposition}\label{bending} {\em There exist positive constants $C_1$ and $\varepsilon_0$, depending only on $M$,  so that for all $k\in\N$ and $0< \varepsilon\leq \varepsilon_0$ we can find a Lipschitz map $F:M\rightarrow M$  such that
\begin{itemize}
\item $F$  is homotopic to the identity;
\item $F(M\setminus B_{\varepsilon3^{-k}}(c(k)))\subset G(K(k)_n)$;
\item $|DF|\leq C_1\varepsilon^{-1}$.
\end{itemize}
}

\begin{proof} Let $x_0$ be the center of the unit cube $I^{n+1}$, and let $\delta$ be a positive constant, to be chosen later. We start by constructing  $f_{\delta}:I^{n+1}\rightarrow I^{n+1}$ a Lipschitz map such that
\begin{itemize}
\item $f_\delta(x)=x$ for every $x\in \partial I^{n+1}\cup\{x_0\}$;
\item $f_\delta$  is homotopic to the identity relative to $\partial I^{n+1}$;
\item $f_\delta(I^{n+1}\setminus B_{\varepsilon}(0))\subset \partial I^{n+1}$;
\item $|Df_\delta|\leq c\delta^{-1}$, where $c=c(n)$.
\end{itemize}
Choose $C$  a bilipschitz homeomorphism between the cube and the unit  ball that sends $x_0$ to the origin. Let $\eta:\mathbb{R}\rightarrow \mathbb{R}$ be a smooth function such that $\eta(t)=1$ if $t\leq 1/2$, $\eta(t)=0$ if $t\geq 1$ and $0\leq \eta(t)\leq 1$ for every $t \in \mathbb{R}$. Set $\eta_\delta(t)=\eta(t/\delta)$ and  
$$h_\delta(x) = \eta_\delta(|x|)x+(1-\eta_\delta(|x|))\frac{x}{|x|},\quad \mbox{for }x \in \overline{B}_1(0).$$ 

The map $f_{\delta}=C^{-1}\circ h_{\delta}\circ C$ satisfies all the required properties.

 For each $\sigma\in K(k)_{n+1}$, we pick an affine linear homomorphism  $L_{\sigma}:I^{n+1}\rightarrow \sigma$ with $L_\sigma(x_0)=q_\sigma$, where $q_\sigma\in I^{n+1}$ denotes the center
 of $\sigma$, and  define  $$F_{\sigma}:G(\sigma)\rightarrow G(\sigma), \quad F_{\sigma}=G\circ L_{\sigma}\circ f_{\delta}\circ L_{\sigma}^{-1}\circ G^{-1}.$$
 
 The map $F_{\sigma}$ satisfies the following conditions:
 \begin{itemize}
 \item $F_\sigma(x)=x$ for every $x\in \partial G(\sigma)$;
\item $F_\sigma$  is homotopic to the identity relative to $\partial G(\sigma)$;
\item $F(G(\sigma)\setminus B_{\delta3^{-k}L^{-1}}(q_\sigma))\subset \partial G(\sigma)$;
\item $|DF|\leq c_{1,\sigma}\delta^{-1}$,
\end{itemize}
where $c_{1,\sigma}>0$ depends only on ${M}$ and $L$ is the Lipschitz constant of $G^{-1}:M \rightarrow K$.
 
We choose $\delta=\varepsilon\, L$, and define $F:M\rightarrow M$ by $F(x)=F_\sigma(x)$ if $x\in \sigma$. The map $F$ is well-defined and satisfies the desired properties.

 \end{proof}

\begin{proof}[Proof of Theorem \ref{upper.bound}] 
Let $p\in\N$. Choose $k\in \N\cup\{0\}$   so that  $3^k\leq p^{\frac{1}{n+1}}\leq 3^{k+1}$.

Let $f: M \rightarrow \mathbb{R}$ be a function satisfying properties (i), (ii) and (iii) of Lemma \ref{morse}. By Lemma \ref{morse.flat} (i), the open set  $\{x\in M: f(x)<t\}$ has finite perimeter for all $t$.
Hence, by \cite[Theorem 30.3]{simon}, we have a well-defined element $$f^{-1}(t)=\partial \{x\in M: f(x)<t\}\in \mathcal Z_n(M;\Z_2).$$ 

For each  $a=(a_0,\ldots,a_p) \in \mathbb{R}^{p+1}$, $|a|=1$, we consider the polynomial $P_a(t)=\sum_{i=0}^p a_i t^i$. Let $t_a^{(1)}, \dots, t_a^{(k_a)}$ be the zeros of $P_a$, where   $k_a \leq p$.

We then  define a  function 
$$\hat{\Psi}:\{a\in \R^{p+1}:|a|=1\}\rightarrow \mathcal{Z}_n(M;\Z_2)$$
by
 $$
\hat{\Psi}(a_0,\ldots,a_p)=\partial\left \{x\in M:P_a(f(x))<0\right\}.$$
Note that the open set $ \{x\in M:P_a(f(x))<0\}$ has finite perimeter, since
\begin{equation}\label{inclusion.perimeter}
\{x\in M:P_a(f(x))=0\}\subset  f^{-1}(t_a^{(1)}) \cup \cdots \cup f^{-1}(t_a^{(k_a)}).
\end{equation}

The fact that we are using  $\Z_2$ coefficients implies that  $\Psi(a)=\Psi(-a)$, and therefore $\hat{\Psi}$ induces a  map $\Psi:\RP^p\rightarrow \mathcal{Z}_n(M;\Z_2).$

\subsection{Claim}{\em  The function $\Psi$ is continuous in the flat topology.}

\medskip
 Let $\{\theta_j\}_{j\in \N}$ be a sequence in $S^p$ that converges to  $\theta\in S^p$. It suffices to show that
 $$
 \lim_{j \rightarrow \infty} {\bf M}\left(\{x\in M:P_\theta(f(x))<0\} \bigtriangleup
\{x\in M:P_{\theta_j}(f(x))<0\} \right)=0,
 $$
 where $X \bigtriangleup Y= (X \setminus Y) \cup (Y\setminus X)$ denotes the symmetric difference of the sets $X$ and $Y$.
 
Since $P_{\theta_j}\circ f$ converges uniformly to $P_\theta\circ f$, it follows that for any $\alpha>0$  we have
 \begin{multline*}
\{x\in M:P_\theta(f(x))<0\} \bigtriangleup
\{x\in M:P_{\theta_j}(f(x))<0\}  \\
\subset  \{x\in M:-\alpha \leq P_\theta(f(x))\leq \alpha\} 
= f^{-1}\left(\{t: P_\theta(t) \in [-\alpha, \alpha] \}\right)
 \end{multline*}
 for all sufficiently large $j$. But
 $$
  \lim_{\alpha \rightarrow 0} {\bf M} \left(f^{-1}\left(P_\theta^{-1}( [-\alpha, \alpha])\right)\right)=0,
 $$
 by item (iii) of Lemma \ref{morse.flat}. This finishes the proof of the claim.
 
 \subsection{Claim}{\em  The function $\Psi$ belongs to $\mathcal P_p$.}
 
 \medskip

The curve $$\gamma:S^1\rightarrow \RP^p,\quad e^{i\theta}\mapsto [(\cos (\theta/2), \sin (\theta/2), 0,\ldots, 0)],$$ is a generator of $\pi_1(\RP^p).$ Then 
$$\Psi\circ\gamma:S^1\rightarrow  \mathcal{Z}_n(M;\Z_2), \quad e^{i\theta}\mapsto\partial\{x\in M:f(x)<-\cot(\theta/2)\},$$
is a  sweepout  of $M$. The generator $\lambda\in H^1(\RP^p;\Z_2)$ satisfies $\lambda(\gamma)=1$ and $\lambda^p\neq 0$, and so $\Psi$ is a $p$-sweepout. Finally, we see from item (ii) of Lemma \ref{morse.flat} and  inclusion \eqref{inclusion.perimeter} that $\Psi$ has no concentration of mass. This finishes the proof that    $\Psi\in \mathcal{P}_p$.

\medskip

 By Lemma \ref{morse} (iii), no point in $c(k)$ is critical for $f$. Hence, if  $\varepsilon$ is chosen sufficiently small we have  that
$${\bf M}(f^{-1}(t)\llcorner B_{\varepsilon3^{-k}}(x))\leq 2\omega_n\varepsilon^n3^{-nk}\quad\mbox{for all }x\in c(k)\mbox{ and }t \in \mathbb{R},$$
where $\omega_n$ is the volume of the unit $n$-ball. By Lemma \ref{morse} (ii), we can also arrange (by choosing $\varepsilon$ even smaller if necessary) that 
$$
f(B_{\varepsilon3^{-k}}(x)) \cap f(B_{\varepsilon3^{-k}}(y)) = \emptyset
$$
for all $x,y \in c(k)$ with $x \neq y$. In particular,
$$
{\bf M}\left (f^{-1}(t) \llcorner
 B_{\varepsilon3^{-k}}(c(k))\right) \leq 2 \omega_n \varepsilon^n3^{-nk}
$$
for every $t \in \mathbb{R}$.

For that choice of $\varepsilon$, we take  the map $F$ given by Proposition \ref{bending} and set $$\Phi:\RP^p\rightarrow  \mathcal{Z}_n(M;\Z_2),\quad \Phi(\theta)=F_{\#}(\Psi(\theta)).$$ Since $F$ is Lipschitz and homotopic to the identity we obtain that $\Phi \in \mathcal P_p$.

We now estimate ${\bf M}(\Phi(\theta))$ for all $\theta\in \RP^p$. We have
\begin{multline*}
{\bf M}\left(F_{\#}\left (f^{-1}(t)\llcorner  B_{\varepsilon3^{-k}}(c(k))\right)\right)\leq (\sup_M |DF|)^n {\bf M}\left(f^{-1}(t)\llcorner  B_{\varepsilon 3^{-k}}(c(k))\right)\\
\leq 2(\sup_M |DF|)^n \omega_n \varepsilon^n3^{-nk}\leq 2C_1^n\omega_n3^{-nk}.
\end{multline*}
Because each $\Psi(\theta)$ consists of at most $p$ level surfaces of $f$, we obtain
\begin{equation}\label{bound1}
{\bf M}\left(F_{\#}\left (\Psi(\theta)\llcorner  B_{\varepsilon3^{-k}}(c(k))\right)\right)
 \leq 2pC_1^n \omega_n3^{-nk}
 \end{equation}
 for all $\theta\in \RP^p$.
 
Set $B=M\setminus  B_{\varepsilon3^{-k}}(c(k))$. From the first property of Proposition \ref{bending} we have that the support of $F_{\#} (\Psi(\theta)\llcorner B )$ is contained in the $n$-skeleton $G(K(k)_n)$. Since  we are using $\Z_2$ coefficients  the multiplicity is at most one. Hence 
$${\bf M}\big(F_{\#} (\Psi(\theta)\llcorner B )\big) \leq {\bf M}(G(K(k)_n))\leq  C_2(\sup_{K}|DG|)^n3^{k(n+1)}3^{-kn}=C_33^{k},$$
where  $C_2$ is  the number of $(n+1)$-cells  in the cell complex $K$ and $C_3=C_2(\sup_{K}|DG|)^n$ depends only on $M$. 

Combining this inequality with \eqref{bound1},  and since $3^k\leq p^{\frac{1}{n+1}}\leq 3^{k+1}$, we have, for some constant $C=C(M)$, 
$${\bf M}(\Phi(\theta))\leq 2pC_1^n\omega_n3^{-nk} + C_33^{k}\leq Cp^{\frac {1}{n+1}}\mbox{ for all }\theta\in \RP^p.$$ 

Therefore $\omega_p(M)\leq C p^{\frac{1}{n+1}}$.
\end{proof}

\section{Equality case}\label{equality}
We apply Lusternik-Schnirelmann theory to prove:

\subsection{Theorem}\label{no.equal}\textit{Assume that $2\leq n\leq 6$ . If $\omega_{p}(M)=\omega_{p+1}(M)$ for some $p\in\N$, then
there exist infinitely many distinct smooth, closed, embedded minimal hypersurfaces in $M$.
}

\begin{proof}
By Proposition \ref{minmax.cohomology}, we can assume that there exist   a (p+1)-admissible cubical subcomplex $X$ and  {a class of $(p+1)$-sweepouts} $$\Pi\in [X,\mathcal{Z}_n(M;{\bf M};\mathbb{Z}_2)]^{\#}$$  so that $\omega_{p+1}(M)={\bf L}(\Pi)$. According to  Proposition \ref{pulltight}, we can  find a critical sequence $S=\{\phi_i\}_{i\in \N}\in \Pi$ so that every $\Sigma \in {\bf C}(S)$ is a stationary varifold with mass equal to ${\bf L}(S)={\bf L}(\Pi)=\omega_{p+1}(M)$.  If $\Phi_i:X\rightarrow \mathcal{Z}_n(M;{\bf M};\Z_2)$ denotes the Almgren extension of $\phi_i$, the fact that $\Pi$ is a {class of $(p+1)$-sweepouts} means that $\Phi_i\in \mathcal{P}_{p+1}$ for all $i$ sufficiently large.

Suppose, by contradiction, that there are only finitely many smooth, closed, embedded minimal hypersurfaces in $M$. Let  $\mathcal S$ be the set of all stationary integral varifolds with area bounded above by $w_{p+1}(M)$ and whose support is a smooth closed embedded hypersurface. 
We consider also the set $\mathcal{T}$ of all mod 2 flat chains $T \in \mathcal{Z}_n(M,\Z_2)$ with ${\bf M}(T) \leq w_{p+1}(M)$ and such that either $T=0$ or the  support of $T$  is a smooth closed embedded minimal hypersurface. By the contradiction hypothesis, both sets $\mathcal{S}$ and $\mathcal{T}$ are finite.

\subsection{Claim}\label{flat.vs.F}\textit{ For every $\varepsilon>0$, there exists $\eta_1>0$ such that 
$$
T \in \mathcal{Z}_n(M,\Z_2)\mbox{ with } {\bf F}(|T|,\mathcal{S}) \leq 2\eta_1 \Longrightarrow  \mathcal{F}(T,\mathcal{T}) < \varepsilon.
$$
}

\begin{proof}
Suppose the claim is false. Then we can find a sequence $\{T_k\} \subset  \mathcal{Z}_n(M,\Z_2)$ with 
${\bf F}(|T_k|,\mathcal{S}) <1/k$  and $\mathcal{F}(T_k,\mathcal{T}) \geq \varepsilon$ for every $k$. By compactness, there exists a subsequence $\{T_l\} \subset \{T_k\}$ that converges in the flat topology to
some $T \in  \mathcal{Z}_n(M,\Z_2)$ and whose associated sequence of varifolds $\{|T_l|\}$ converges in varifold 
topology to some $V \in \mathcal{S}$. In particular, $\mathcal{F}(T,\mathcal{T}) \geq \varepsilon$ and ${\bf M}(T)\leq \omega_{p+1}(M)$.  We also have,
by lower semicontinuity of mass,  that $${\bf M}\left( T \llcorner (M\setminus\mathrm{spt}||V||) \right)  =0.$$

This implies that the support of $T$ is contained in the smooth, closed, embedded minimal hypersurface $\mathrm{spt}||V||$. By the Constancy Theorem (\cite{simon}), $T \in \mathcal{T}$. This is a contradiction, since
$\mathcal{F}(T,\mathcal{T}) \geq \varepsilon$.

\end{proof}

By Proposition \ref{locally.trivial}, there exists $\varepsilon>0$ such that every map $\Phi:S^1\rightarrow \mathcal{Z}_n(M;\Z_2)$ with
$$\Phi(S^1) \subset B_\varepsilon^\mathcal{F}(\mathcal{T}) = \{T \in \mathcal{Z}_n(M;\Z_2): \mathcal{F}(T,\mathcal{T})< \varepsilon\}$$ is trivial. For this given $\varepsilon$, we choose $\eta_1$ as in Claim \ref{flat.vs.F}.

With $k_i\in\N$ so that $\mbox{dmn}(\phi_i)=X(k_i)_0$, consider $Y_i$ to be the cubical subcomplex of $X(k_i)$ consisting of all cells $\alpha \in X(k_i)$ so that 
$${\bf F}(|\phi_i(x)|,\mathcal S)\geq \eta_1$$
for every vertex $x$  in $\alpha_0$. In particular $Y_i$ is a cubical subcomplex of $I(m,k_i)$ for some $m\in\N$. It also follows that 
 \begin{equation}\label{eta}
 {\bf F}(|\Phi_i(x)|,\mathcal S)< 2\eta_1\mbox{ for every }x\in X\setminus Y_i
 \end{equation} if $i$ is sufficiently large.

\subsection{Claim}\label{detect.p-class}\textit{For all  $i$ sufficiently large we have
 $(\Phi_i)_{|Y_i} \in \mathcal{P}_p$. 
}

\begin{proof}
Assume $i$ is sufficiently large so that $\Phi_i\in \mathcal{P}_{p+1}$ and \eqref{eta} holds.

The map $(\Phi_i)_{|Y_i}$ is continuous in the flat topology and has no concentration of mass (Lemma \ref{mass.concentration.lemm}) and thus we only need to check that it is a $p$-sweepout.

Let $\lambda=\Phi_i^*(\overline{\lambda}) \in H^1(X;\Z_2)$. Then,  since  $\Phi_i$ is a  $(p+1)$-sweepout (see Definition \ref{defi.pclass}), we have
\begin{itemize}
\item  for every curve $\gamma:S^1 \rightarrow X$ we have $\lambda(\gamma)\neq 0$ if and only if $\Phi_i\circ\gamma$ is a sweepout;
\item $\lambda^{p+1} \neq 0$ in  $H^{p+1}(X;\Z_2)$.
\end{itemize}

 Let  $Z_i=\overline{X\setminus Y_i}$.  Hence $Z_i$ is a subcomplex of $X(k_i)$ as well. Consider the inclusion maps $i_{1}: Z_i\rightarrow X$ and $i_{2}: Y_i\rightarrow X$. 
 
 If we show that $(i^{*}_{2}\lambda)^p \neq 0$ in $H^p(Y_i;\Z_2)$, it follows at once that $(\Phi_i)_{|Y_i}$ is a $p$-sweepout.

	For any closed curve $\gamma:S^1 \rightarrow Z_i$, we have  from Claim \ref{flat.vs.F} and \eqref{eta} that
	$$\Phi_i\circ\gamma (S^1) \subset B_\varepsilon^\mathcal{F}(\mathcal{T}).$$
Proposition \ref{locally.trivial}  implies that $\Phi_i\circ\gamma:S^1 \rightarrow \mathcal{Z}_n(M;\Z_2)$ is trivial and, as  a result, $i^*_{1}\lambda(\gamma)=0$. This means    $i_{1}^*\lambda=0$ in $H^1(Z_i;\Z_2)$ because $H^1(Z_i;\Z_2)={\rm Hom}\,(H_1(Z_i);\Z_2)$, by the Universal Coefficient Theorem.

 From  the natural exact sequence
	$$H^1(X,Z_i;\Z_2)\mathop{\rightarrow}^{j^*} H^1(X;\Z_2) \mathop{\rightarrow}^{i_{1}^*} H^1(Z_i ;\Z_2)$$
we obtain  that $\lambda=j^{*}\lambda_1$ for some $\lambda_1\in H^1(X,Z_i;\Z_2)$.

	Suppose  $i_2^{*}(\lambda^p)=0$. Then the exact sequence
		$$H^p(X,Y_i;\Z_2)\mathop{\rightarrow}^{j^*} H^p(X;\Z_2) \mathop{\rightarrow}^{i_2^*} H^p(Y_i ;\Z_2)$$
implies that $j^*\lambda_2=\lambda^p$ for some $\lambda_2\in H^p(X,Y_i;\Z_2)$. 

Thus
$$ j^{*}\lambda_1 \smile j^*\lambda_2=\lambda^{p+1}\neq 0\mbox{ in }H^{p+1}(X;\Z_2).$$

On the other hand, since $Y_i$ and $Z_i$ are subcomplexes of $X(k_i)$, there is a natural notion of relative cup product (see \cite{hatcher}, p 209):
$$
H^1(X,Z_i;\Z_2) \smile H^p(X,Y_i;\Z_2) \rightarrow H^{p+1}(X,Y_i\cup Z_i ;\Z_2).
$$
But $Y_i\cup Z_i =X$, hence $H^{p+1}(X,Y_i\cup Z_i ;\Z_2)=H^{p+1}(X,X ;\Z_2)=0$. In particular, $\lambda_1 \smile \lambda_2=0$. This is a contradiction because $$j^*(\lambda_1 \smile \lambda_2) =  j^{*}\lambda_1 \smile j^*\lambda_2=\lambda^{p+1}\neq 0.$$ Hence $i_{2}^{*}(\lambda^p)\neq0$ and the proof is finished.


\end{proof}

Consider the sequence $\tilde{S}=\{\psi_i\}$, where 
$$\psi_i=(\phi_i)_{|Y_i}:(Y_i)_0 \rightarrow \mathcal{Z}_n(M;\Z_2),$$
and let $$L={\bf L}(\tilde{S})=\limsup_{i\to\infty}\max\{{\bf M}(\psi_i(y)):y\in (Y_i)_0\}.$$ Of course $L\leq \omega_{p+1}(M).$ There are two cases to consider: $L<\omega_{p+1}(M)$ and $L=\omega_{p+1}(M)$. 

If  $L<\omega_{p+1}(M)$, then by property (iii) of  Theorem \ref{discrete.continuous} we have that the Almgren extension $\Phi_i$ satisfies 
$$\sup_{y\in Y_i} {\bf M}(\Phi_i(y))<\omega_{p+1}(M)$$ for sufficiently large $i$. On the other hand, we know from Claim \ref{detect.p-class}  that $(\Phi_i)_{|Y_i}\in\mathcal P_p$ and thus 
$$\sup_{y\in Y_i} {\bf M}(\Phi_i(y))\geq \omega_p(M)=\omega_{p+1}(M),$$
which is a contradiction.

Suppose now that $L=\omega_{p+1}(M)$. Since
\begin{multline*}
{\bf C}(\tilde{S})=\{V:||V||(M)=L,\, V=\lim_{j\to\infty}|\psi_{i_j}(y_j)|\mbox{ as varifolds,} \\
\mbox{ for some increasing sequence }\{i_j\}_{j\in \N} \mbox{ and }y_j\in \mathrm{dmn}(\psi_{i_j})\},
\end{multline*}
we have that  ${\bf C}(\tilde{S}) \subset {\bf C}(S)$. We also have that 
${\bf C}(\tilde{S}) \subset \{V: {\bf F}(V, \mathcal{S}) \geq \eta_1\}$, by definition of $Y_i$. We conclude that
although every element of ${\bf C}(\tilde{S})$ is stationary, none of them has smooth support. In particular no element of ${\bf C}(\tilde{S})$ is $\mathbb{Z}_2$ almost minimizing in annuli. 

Therefore we can apply Theorem \ref{combinatorial.pitts} and
produce a sequence $\tilde{S}^*=\{\psi_i^*\}$ of maps 
$$\psi_i^*: Y_i(l_i)_0 \rightarrow \mathcal{Z}_n(M;\Z_2)$$ such that:
\begin{itemize}
\item $\psi_i$ and $\psi_i^*$ are homotopic to each other with finesses that tend to zero as $i \rightarrow \infty$,
\item ${\bf L}(\tilde{S}^*)=\limsup_{i\to\infty}\max\{{\bf M}(\psi_i^*(y)):y\in Y_i(l_i)_0\}<{\bf L}(\tilde{S})=L.$
\end{itemize}

By Proposition \ref{discrete.homotopy.continuous}, the first item above implies that the Almgren extensions $\Psi_i$, $\Psi_i^*$ to $Y_i$ of $\psi_i$, $\psi_i^*$, respectively, are homotopic to each other if $i$ is sufficiently large.  Moreover, Theorem \ref{discrete.continuous} (ii) implies that $\Psi_i=(\Phi_i)_{|Y_i}$ and thus we have from Claim  \ref{detect.p-class} that $\Psi_i^*\in\mathcal P_p$ for all $i$ sufficiently large. Hence
$$\sup_{y\in Y_i} {\bf M}(\Psi_i^*(y))\geq \omega_p(M)$$
for all large $i$. The second item implies, by property (iii) of  Theorem \ref{discrete.continuous}, that $$\sup_{y\in Y_i} {\bf M}(\Psi_i^*(y))<L=\omega_{p+1}(M)=\omega_p(M)$$ for all large $i$ and we get a contradiction.

Both  cases $L<\omega_{p+1}(M)$ and $L=\omega_{p+1}(M)$ lead to  a contradiction, hence there must be infinitely many distinct smooth, closed, embedded minimal hypersurfaces in $M$.

 \end{proof}

\section{Proof of Main Theorem}\label{proof.main.theorem}

By contradiction, suppose that the set $\mathcal L$ of all smooth, connected, closed, embedded minimal hypersurfaces of $M$ is finite and that any disjoint subcollection of $\mathcal L$ has at most $n$ elements.

It follows from Proposition \ref{minmax.cohomology} that for every $p\geq 1$ we can  find  $p$-admissible cubical subcomplexes $X_p$ and   $\Pi_p\in  [X_p,\mathcal{Z}_n(M;{\bf M};\mathbb{Z}_2)]^{\#}$ so that 
$$\omega_p(M)={\bf L}(\Pi_p).$$  By Theorem \ref{pitts.min.max} and Theorem \ref{regularity.thm}, we have  
$$\omega_p(M)=||V_p||(M)$$
for some $V_p \in \mathcal{IV}_n(M)$, where $V_p$ is the varifold of a smooth, closed, embedded minimal hypersurface, with possible multiplicities.

We can write
 $$
 V_p= n^{(p)}_1 \Sigma^{(p)}_1 + \cdots + n^{(p)}_{l_p} \Sigma^{(p)}_{l_p}
 $$
 with $\Sigma^{(p)}_j \in \mathcal L$, $n^{(p)}_j \in \mathbb{N}$ for $1\leq j\leq l_p$.  Since the support of $V_p$ is embedded, we have $\{\Sigma^{(p)}_1, \dots, \Sigma^{(p)}_{l_p}\}$ is disjoint and hence $l_p \leq n$.
 
Because we are assuming that $\mathcal{L}$ is finite, we must have by  Theorem \ref{no.equal}  that $$\omega_p=||V_p||(M)<||V_{p+1}||(M)=\omega_{p+1}\quad\mbox{for all }p\in\N.$$
Hence 
$$
\#\{\omega_k(M):k=1,\ldots, p\}= p.
$$

Let $\delta>0$ be such that $|\Sigma| \geq \delta$ for every $\Sigma \in \mathcal{L}$. By Theorem \ref{upper.bound} one has $\omega_p(M) \leq Cp^\frac{1}{n+1}$, and then 
$n^{(p)}_j \in \{1,\dots,\lfloor Cp^\frac{1}{n+1}/\delta\rfloor\}$. This implies 
$$
\#\{\omega_k(M):k=1,\ldots, p\}\leq C'p^\frac{n}{n+1}
$$
for a constant $C'>0$ independent of $p$. We get a contradiction when $p$ is large, and this finishes the proof.

\section{Lower Bounds}\label{lower.bounds}

The following result was proven by Gromov (see \cite[Section 4.2.B]{gromov0} or \cite[Section 8]{gromov}). For the convenience of the reader we present a proof of this theorem that follows closely the proof  given by Guth in \cite[Section 3]{guth}.

\subsection{Theorem}\label{gromov}{\em There exists  $C=C(M)>0$ so that 
$$\omega_p(M)\geq C p^{\frac {1}{n+1}}\quad\mbox{for all }p\in\N.$$
}


Given $p\in M$, let $B_r(p)$ denote the geodesic ball in $M$ of radius $r$ and centered at $p$.

\subsection{Proposition}\label{lowerbound.lemma}{\em There exist  positive constants $\alpha_0=\alpha_0(M)$ and $r_0=r_0(M)$ so that for any sweepout $\Phi:S^1\rightarrow \mathcal{Z}_n(M;\Z_2)$,  we have
$$\sup_{\theta \in S^1}{\bf M}( \Phi(\theta)\llcorner B_r(x))\geq \alpha_0r^n$$
for all $x\in M$ and $0< r\leq r_0$.
}
\begin{proof}
We will use notation and definitions of Section \ref{almgren.homo.discrete}.

The compactness of $M$ and scaling considerations imply we can find positive constants $\rho_1$ and $r_1$, depending only on $M$,  so that
$$\nu_{B_r(x),\partial B_r(x)}>\alpha_1r^n\quad\mbox{for all}\quad x\in M\mbox{ and } 0<r\leq r_1,$$
This means that for all $$T  \in \mathcal{Z}_n(B_{r}(x), \partial B_{r}(x);\Z_2)\quad\mbox{with}\quad
\mathcal{F}(T)<\alpha_1r^{n+1},$$
  there exists an isoperimetric choice $Q\in I_{n+1}(B_r(x);\Z_2)$ with
$$
\partial Q - T \in I_n(\partial B_r(x);\Z_2),
$$
that is unique assuming ${\bf M}(Q)<\alpha_1r^{n+1}$ (Lemma \ref{constancy.lemma}).

Let  $x\in M$ and $0<r\leq r_1$. Choose $\delta$ small so that $(1+\frac{2}{r}) \rho\delta< \alpha_1\left(\frac{r}{2}\right)^{n+1}$ and $k$ sufficiently large so that
$$\mathcal{F}(\Phi(e^{2\pi ix}),\Phi(e^{2\pi iy}))\leq \delta\quad\mbox{for all }x, y \mbox{ in some common cell of  }I(1,k).
$$
We set
 $$\phi:I(1,k)_0\rightarrow \mathcal{Z}_n(M;\Z_2)\quad\phi([x])=\Phi(e^{2\pi ix}).$$ 
Assuming $\delta<\nu_{M}$,  we can find an isoperimetric choice $Q_j \in {\bf I}_{n+1}(M;\Z_2)$, $j=0,\ldots,3^{k}-1$, such that 
$$\partial Q_j=\phi(a_{j+1})-\phi(a_j)\quad \mbox{and}\quad {\bf M}(Q_j)\leq \
\rho{\mathcal F}(\phi(a_{j+1})-\phi(a_j))\leq \rho\delta,$$
where   $a_j=[j3^{-k}]$ and $\rho=\rho(M)$ is defined in Section \ref{almgren.homo.discrete}.
The fact that $\Phi$ is a  sweepout implies that we can also assume that
\begin{equation}\label{fund.cycle}
 \sum_{j=0}^{3^k-1}Q_j= M\quad\mbox{in }{\bf I}_{n+1}(M;\Z_2).
\end{equation}

We can find $r/2\leq s\leq r$ (\cite[Lemma 28.5]{simon}) so that
$$\phi(a_{j})\llcorner B_s(x)\in \mathcal{Z}_n(B_s(x),\partial B_s(x);\Z_2),$$ 
$$L_j=\partial \left(Q_j\llcorner B_s(x)\right)-\partial Q_j \llcorner B_s(x) \in {\bf I}_n(\partial B_s(x);\Z_2),$$
and such that
$${\bf M}(L_j)\leq \frac{2}{r}{\bf M}(Q_j)\quad\mbox{for all}\quad j=0,\ldots,3^k-1.$$
Let $$\bar \phi:I(1,k)_0\rightarrow  \mathcal{Z}_n(B_s(x),\partial B_s(x);\Z_2),\quad \bar \phi(x)=\phi(x)\llcorner B_s(x).$$
  Since
 \begin{multline*}
\mathcal F(\bar \phi(a_{j+1})-\bar \phi(a_{j})) \leq {\bf M} (Q_j\llcorner B_s(x))+{\bf M}(L_j)\leq \left(1+\frac 2 r\right){\bf M} (Q_j)\\
\leq  \left(1+\frac 2 r\right) \rho\delta<\alpha_1 \left(\frac{r}{2}\right)^{n+1}<\alpha_1s^{n+1}
\end{multline*}
and $${\bf M}(Q_j\llcorner B_s(x))\leq {\bf M}(Q_j)\leq\rho\delta<\alpha_1s^{n+1},$$
we have that $Q_j\llcorner B_s(x)$ is the isoperimetric choice for $\phi(a_{j+1})-\bar \phi(a_{j})$.
Therefore, recalling the definition in \ref{almgren.homo.boundary},
\begin{equation}\label{trivial}
 F_{B_s(x),\partial B_s(x)}^{\#}(\bar\phi)=\left[\sum_{j=0}^{3^k-1} Q_j\llcorner B_s(x)\right]=[M\llcorner B_s(x)]=[B_s(x)].
 \end{equation}

From \cite[Proposition 1.22]{almgren}, using the compactness of $M$ and scaling considerations,  we can choose $\alpha_2>0$ and $\rho_2>0$ depending only on $M$  so that for each $x\in M$, $0<r\leq r_1$ and 
$$T  \in \mathcal{Z}_n(B_{r}(x), \partial B_{r}(x);\Z_2)\quad\mbox{with}\quad{\bf M}(T)<\alpha_2 r^{n},$$
 there exists $Q\in I_{n+1}(B_r(x);\Z_2)$ with
$$
\partial Q - T \in I_n(\partial B_r(x);\Z_2)  \quad\mbox{and}\quad {\bf M}(Q) \leq \rho_2{\bf M}(T)^\frac{n+1}{n}.
$$
Set $\alpha_0=\min\{\alpha_2,\alpha_1/(2\rho_2)\}.$

\medskip
\noindent{\bf Claim:} There exists $x\in I(1,k)_0$ such that ${\bf M}(\bar \phi(x))\geq \alpha_0 {s}^{n}$.

\medskip

Suppose, by contradiction, that the claim is false. Then ${\bf M}(\bar \phi(x))< \alpha_0 s^{n}$ for all $x\in I(1,k)_0$.
This implies we can find $S_j\in {\bf I}_{n+1}( B_s(x);\Z_2)$, for all  $j=0,\ldots,3^k$, so that
$$\partial S_j-\bar \phi(a_{j})\in  {\bf I}_{n}(\partial  B_s(x);\Z_2)\quad\mbox{and}\quad {\bf M}(S_j)< \rho_2 {\bf M}(\bar \phi(a_j))^\frac{n+1}{n}.$$
Note that $S_{3^k}=S_0$ because $\bar\phi([0])=\bar\phi([1])$. 

Furthermore, $S_{j+1}-S_{j}$ is also an isoperimetric choice for $\phi(a_{j+1})-\bar \phi(a_{j})$.  It must be  equal to $Q_j\llcorner B_s(x)$ because
$${\bf M}(S_{j+1}-S_{j})\leq \rho_2 {\bf M}(\bar \phi(a_{j+1}))^\frac{n+1}{n}+\rho_2 {\bf M}(\bar \phi(a_j))^\frac{n+1}{n}<2\rho_2\alpha_0s^{n+1}\leq \alpha_1s^n. $$
As  a result,
$$ F_{B_s(x),\partial B_s(x)}^{\#}(\bar\phi)=\left[\sum_{j=0}^{3^k-1} Q_j\llcorner B_s(x)\right]=S_{3^k}-S_0=0.$$
This contradicts \eqref{trivial} and thus proving the claim.

\medskip

The claim implies the existence of some $\theta\in S^1$ with 
$${\bf M}( \Phi(\theta)\llcorner B_r(x))\geq 2^{-n}\alpha_0r^n.$$

\end{proof}

\begin{proof}[Proof of Theorem \ref{gromov}]
By Proposition \ref{discrete.continuous.homotopy} (ii), it suffices to show that for every $p$-admissible $X$ and every $p$-sweepout $\Phi:X\rightarrow \mathcal{Z}_n(M;{\bf M};\Z_2)$  continuous in the mass topology, we have
$$\sup_{x\in X}\,{\bf M}(\Phi(x))\geq Cp^{\frac{1}{n+1}},$$
where $C$ is a positive constant that depends only on $M$.

There exists some  constant $\nu=\nu(M)>0$ such that, for every $p\in \N$, one  can find  a collection of $p$ disjoint geodesic balls $\{B_j\}_{j=1}^p$ of radius $r=\nu p^{-\frac{1}{n+1}}$. 
Let $\alpha_0>0$ be  the constant of Proposition \ref{lowerbound.lemma}.

Fix $p\in \mathbb{N}$. We can choose $k$ sufficiently large so that $${\bf M}(\Phi(x),\Phi(y))<\frac{\alpha_0}{6}r^n$$
for all $x,y$ in some common cell of $X(k)$.
We define $S_j$ as the union of all cells $\sigma$ of $X(k)$ so that 
$$
{\bf M}(\Phi(x)\llcorner B_j)\leq \frac{\alpha_0}{3} r^n
$$
for every $x\in \sigma_0$. In particular, ${\bf M}(\Phi(y)\llcorner B_j)<\frac{\alpha_0}{2} r^n$ for every $y\in S_j$.

\subsection{Lemma}\label{gromov.lowerbound}\textit{There exists $x\in X\setminus (S_1\cup\ldots\cup S_p)$}.
\begin{proof}
Suppose  $X=S_1\cup\ldots\cup S_p$, by contradiction.  

Since  $\Phi$ is a  $p$-sweepout we have, with $\lambda=\Phi^*(\overline{\lambda}) \in H^1(X;\Z_2)$, that 
\begin{itemize}
\item  for every curve $\gamma:S^1 \rightarrow X$, $\lambda(\gamma)\neq 0$ if and only if $\Phi\circ\gamma$ is a  sweepout;
\item $\lambda^{p} \neq 0$ in  $H^{p}(X;\Z_2)$.
\end{itemize}

 We are going to find a  closed curve $\gamma:S^1 \rightarrow X$ such that $\gamma(S^1)$ is contained in some $S_j$ and so that $\lambda(\gamma)\neq 0$. 
In that case we get that  $\Phi\circ \gamma:S^1\rightarrow \mathcal Z_n(M;{\bf M};\Z_2)$  is a sweepout with ${\bf M}(\Phi(y)\llcorner B_j)<\frac{\alpha_0}{2} r^n$, contradicting Proposition \ref{lowerbound.lemma} applied to the ball $B_j$.

Consider the inclusion maps $i_{S_j}:S_j\rightarrow X$, $j=1, \ldots, p$. 
\subsection{Claim}{\em 
For some $j=1,\ldots,p$, we have $i_{S_j}^{*}(\lambda)\neq 0$ in $H^1(S_j,\Z_2)$.}

\medskip

Suppose $i_{S_j}^{*}(\lambda)= 0$ for all $j=1,\ldots,p$.  Consider the exact sequence
$$H^1(X,S_j;\Z_2)\mathop{\rightarrow}^{j^*} H^1(X;\Z_2) \mathop{\rightarrow}^{i_{S_j}^*} H^1(S_j ;\Z_2).$$
Then we can find $\lambda_j\in H^1(X,S_j;\Z_2)$ so that $j^*(\lambda_j)=\lambda$. Therefore
$$ j^*(\lambda_1)\smile \ldots\smile j^*(\lambda_p)=\lambda^p\neq 0\mbox{ in } H^p(X;\Z_2).$$

Since $S_j$ is a subcomplex of $X(k)$ for each $j$, we have a natural notion of relative cup product (see \cite{hatcher}, p 209):
$$
H^1(X,S_1;\Z_2) \smile \cdots \smile H^1(X,S_p;\Z_2) \rightarrow H^{p}(X,S_1\cup \cdots \cup S_p ;\Z_2).
$$
But we are assuming that $S_1\cup \cdots \cup S_p =X$, hence $$H^{p}(X,S_1\cup \cdots \cup S_p ;\Z_2)=H^p(X,X;\Z_2)=0.$$ 

Therefore 
$$
\lambda^p=j^*(\lambda_1)\smile \ldots\smile j^*(\lambda_p)=j ^*(\lambda_1 \smile \ldots \smile \lambda_p) =0.
$$
This cannot be true, hence  $i_{S_j}^{*}(\lambda)\neq 0$  for some $j=1,\ldots,p$. This proves the claim.

\medskip

Let  $S_j$ be as in the above claim. By the Universal Coefficient Theorem, we have that ${\rm Hom}\,(H_1(S_j);\Z_2)=H^1(S_j;\Z_2)$. Thus we can find a   closed curve $\gamma\subset S_j$ such that  $\lambda(i_{S_j}\circ\gamma)=(i_{S_j}^*\lambda)(\gamma) \neq 0$.  Therefore $i_{S_j}\circ\gamma$ is a  sweepout in $X$, which is exactly what we wanted to prove.
\end{proof}

The lemma we just proved gives the existence of $x\in X\setminus (S_1\cup\ldots\cup S_p)$. Then, from the definition of the sets $S_j$, we get
$${\bf M}(\Phi(x))\geq \sum_{j=1}^p {\bf M}(\Phi(x)\llcorner B_j)\geq p \frac{\alpha_0}{6} r^n\geq \frac {\alpha_0}{6} \nu^n p^{\frac{1}{n+1}}=Cp^{\frac{1}{n+1}},$$
where $C$ is a positive constant that depends only on $M$. This finishes the proof of the theorem.

\end{proof}

\section{Open problems}

{In this section we state and propose some questions regarding min-max theory applied to the class $\mathcal P_p$ of $p$-sweepouts.  

We start by recalling the min-max definition of the $p^{th}$-eigenvalue of $(M,g)$. Set $V= W^{1,2}(M)\setminus\{0\}$ and consider the Rayleigh quotient
$$ E:V\rightarrow [0,\infty],\quad E(f)=\frac{\int_M|\nabla f|^2 dV_g}{\int_M f^2 dV_g}.$$
Then
$$\lambda_p=\mathop{\inf}_{(p+1)-\mbox{plane }P\subset V}\max_{f\in P}E(f).$$
Hence, in light of Definition \ref{p.width}, one can see $\{\omega_p(M)\}_{p\in\N}$ as a nonlinear analogue of the Laplace spectrum of $M$, as proposed by  Gromov \cite{gromov0}. Many interesting problems can be raised out of this analogy.

For instance, Gromov conjectured in \cite[Section 8]{gromov} (also \cite[Section 5.2]{gromov2}) that the sequence $\{\omega_p(M)\}_{p\in\N}$ satisfies a Weyl Law, meaning that \begin{equation}\label{weyl}
\lim_{p\to\infty}\omega_p(M)p^{-\frac {1}{n+1}}=a(n)({\rm vol}(M,g))^{\frac{n}{n+1}},
\end{equation}
where  $a(n)$ is a constant that depends only on $n$. The authors and  Liokumovich  confirmed this conjecture in  \cite{liokumovich-marques-neves}. Note that from Theorem \ref{upper.bound} and Theorem \ref{gromov} we know that the sequence $\{\omega_p(M)p^{-\frac {1}{n+1}}\}_{p\in \N}$ is contained in some compact interval 
$[c_1,c_2] \subset (0,\infty)$.

This analogy can also be  put forward by considering sweepouts whose surfaces are zero sets of linear combinations 
of eigenfunctions. If $\phi_0, \dots, \phi_p$ denote the first $(p+1)$-eigenfunctions for the Laplace operator of $(M,g)$,
where $\phi_0$ is the constant function, we can consider the map 
\begin{align*}
\Phi_p&:\RP^p\rightarrow\mathcal Z_n(M;\Z_2),\\
\Phi_p&([a_0,\ldots,a_p])=\partial\{x\in M:a_0\phi_0(x)+\ldots+a_p\phi_p(x)<0\}.
\end{align*}

It is interesting to compute the numbers $\omega_p(M)$ in specific examples.
For the case of the unit $3$-sphere $S^3$ with the standard metric, we can choose $\phi_1,\phi_2,\phi_3,\phi_4$ to be the coordinate functions and so it is simple to see that
$$\omega_1(S^3)=\omega_2(S^3)=\omega_3(S^3)=\omega_4(S^3)=\max_{\theta \in \RP^4} {\bf M}(\Phi_4(\theta))=4\pi.$$ 
 Note that the Clifford torus is the nodal set of $\phi_5=x_1^2+x_2^2-x_3^2-x_4^2$.
The space of spherical harmonics in $S^3$ of degree less than or equal to $2$ has dimension $14$. For every $\theta\in \RP^{13}$, we have that $\Phi_{13}(\theta)$ intersects almost every closed geodesic in $S^3$ at most $4$ times and so Crofton's formula implies that ${\bf M}(\Phi_{13}(\theta))\leq 8\pi$. Thus
$$\omega_{13}(S^3)\leq \sup_{\theta\in \RP^{13}}{\bf M}(\Phi_{13}(\theta))=8\pi.$$

Nurser \cite{nurser} used the canonical family  found by the authors in \cite{marques-neves} to show that $\omega_5(S^3)=\omega_6(S^3)=\omega_7(S^3)=2\pi^2$ and that $2\pi^2<\omega_9(S^3)<8\pi$. It would be nice to know  for which values of $k$ we have $2\pi^2<\omega_k(S^3)<8\pi$ and whether they are achieved by interesting minimal surfaces. 

The similar problem for $S^2$ seems to be more tractable and Aiex showed in \cite{aiex} that $\omega_i(S^2)=2\pi$ if $i=1,2,3$ and $\omega_i(S^2)=4\pi$ if $i=4,5,6,7,8$. He also computed these widths on some ellipsoids.


Note that a conjecture of Yau \cite{yau1} states that $$c^{-1}\sqrt \lambda_p\leq  \mathcal H^n(\{\phi_p=0\})\leq c\sqrt \lambda_p,$$ 
where $c=c(M,g)>0$. This conjecture was proven by Donnelly and Fefferman \cite{donnelly-fefferman} when the metric is analytic and the lower bound has been recently proved for smooth metrics by Logunov (\cite{logunov}). Note that from Theorem \ref{gromov} one should have
$$\sup_{\theta \in \RP^p}{\bf M}(\Phi_p(\theta))\geq c^{-1}p^{\frac {1}{n+1}}.$$
Assuming  a more speculative nature,  it would be interesting to see if the family $\Phi_p$ defined above is asymptotically optimal.


It is interesting to study the general behavior of the minimal hypersurfaces that are produced by applying min-max theory to the classes $\mathcal P_p$. Is it possible to analyze their Morse indices (see work \cite{marques-neves-index} of the authors)? Do their volumes (not counting multiplicity) become unbounded? How are they distributed? One could naively expect that under generic conditions they should have index $p$, multiplicity one and their volumes converge to infinity. The proof of Theorem \ref{gromov} suggests  that these surfaces might  become equidistributed in space. 
}

\appendix

\section{}\label{A.section}

\begin{proof}[Proof of  Proposition \ref{close.implies.homotopic}.]

It follows from the work of Almgren (\cite{almgren}, Theorem 8.2) that there exist  $0<\delta_0<\ldots<\delta_{m+1}$, depending only on $M$ and $m$, such that if $\Phi:I^k \rightarrow \mathcal{Z}_n(M;\Z_2)$, $k \leq m$, is continuous in the flat topology, $\Phi(x)=0$ for all $x \in \partial I^k$ and $\mathcal{F}(\Phi(x))\leq\delta_k$ for every $x\in I^k$, then there exists a homotopy $H:I^{k+1}\rightarrow \mathcal{Z}_n(M;\Z_2)$ with the following properties:
\begin{itemize}
\item $H$ is continuous in the flat topology;
\item $H(x,0)=0$ and $H(x,1)=\Phi(x)$ for every $x \in I^k$;
\item $H(x,t)=0$ for every $x\in \partial I^k$ and $t\in [0,1]$;
\item $\sup\{\mathcal{F}(H(w)):w\in I^{k+1}\} \leq \delta_{k+1}$.
\end{itemize}

Set $\delta=\delta_0$ and let $\Psi=\Phi_2-\Phi_1$. Denote by $Y^{(j)}$  the union of all cells of $Y$ with dimension at most $j$, respectively, for every $j=0,\dots,m$. We will construct the homotopy by an inductive process.

\subsection{Claim}{\em For each  $j=0,\dots,m$, there exists a map $H:Y^{(j)}\times I \rightarrow \mathcal{Z}_n(M;\Z_2)$ that satisfies:
\begin{itemize}
\item $H$ is continuous in the flat topology;
\item $H(y,0)=0$ and $H(y,1)=\Psi(y)$ for every $y \in Y^{(j)}$;
\item $\sup\{\mathcal{F}(H(w)):w\in Y^{(j)}\times I \} \leq\delta_{j+1}.$
\end{itemize}
}

\medskip

The proof is by induction. Almgren's construction described above gives a map $H:Y^{(0)}\times I \rightarrow \mathcal{Z}_n(M;\Z_2)$ that satisfies 
\begin{itemize}
\item $H$ is continuous in the flat topology;
\item $H(y,0)=0$ and $H(y,1)=\Psi(y)$ for every $y \in Y^{(0)}$;
\item $\sup\{\mathcal{F}(H(w)):w\in Y^{(0)}\times I \} \leq\delta_1.$
\end{itemize}

Let us suppose now that we have constructed a map $H:Y^{(j-1)}\times I \rightarrow \mathcal{Z}_n(M;\Z_2)$ that satisfies 
\begin{itemize}
\item $H$ is continuous in the flat topology;
\item $H(y,0)=0$ and $H(y,1)=\Psi(y)$ for every $y \in Y^{(j-1)}$;
\item $\sup\{\mathcal{F}(H(w)): w\in Y^{(j-1)}\times I \} \leq  \delta_j.$
\end{itemize}
We can extend $H$ continuously to $Y^{(j)} \times \{1\}$ by putting $H(y,1)=\Psi(y)$ for each $y\in Y^{(j)}$, and we will still have $$\sup\{\mathcal{F}(H(w)):w\in (Y^{(j-1)}\times I) \cup  (Y^{(j)} \times \{1\})\} \leq  \delta_j.$$

Let $\sigma\in Y^{(j)}_j$ be a $j$-dimensional cell of $Y$ and choose a homeomorphism $f_\sigma:I^{j+1}\rightarrow \sigma \times I$ such that $f_\sigma(I^j\times \{1\})=(\sigma \times \{1\}) \cup (\partial \sigma \times I).$ Then $H\circ f_\sigma$ is well-defined on $I^j\times \{1\}$. Since $f_\sigma(\partial (I^j\times \{1\}))\subset \partial \sigma \times \{0\}$, then $(H\circ f_\sigma)(x)=0$ for all $x\in \partial (I^j\times \{1\})$. The Almgren's construction  gives again a map $H_\sigma:I^j\times I \rightarrow \mathcal{Z}_n(M;\Z_2)$ that satisfies:
\begin{itemize}
\item $H_\sigma$ is continuous in the flat topology;
\item $H_\sigma(x,0)=0$ and $H_\sigma(x,1)=(H\circ f_\sigma)(x)$ for every $x \in I^j$;
\item $H_\sigma(x,t)=0$ for every $x\in \partial I^j$ and $t\in [0,1]$;
\item $\sup\{\mathcal{F}(H_\sigma(w)):w\in I^j\times I \} \leq  \delta_{j+1}.$
\end{itemize}

We can extend $H$ to a  map $H:Y^{(j)}\times I \rightarrow \mathcal{Z}_n(M;\Z_2)$ by setting 
$H=H_\sigma\circ f_\sigma^{-1}$ on each $\sigma \times I$, $\sigma\in Y^{(j)}_j$. This proves the claim.

\medskip

By applying the claim with $j=m$, we get a homotopy $H$  between the zero map and $\Psi=\Phi_2-\Phi_1$. Then
$\tilde{H}(z)=H(z)+\Phi_1(z)$ for $z \in Y\times I$ is the desired homotopy.
\end{proof}


\bibliographystyle{amsbook}

\end{document}